\documentclass[11pt,reqno]{amsart}
\usepackage[T1]{fontenc}
\usepackage[utf8]{inputenc}
\usepackage{lmodern}
\usepackage[a4paper,margin=30mm]{geometry}
\usepackage{amsmath,amssymb,amsthm,mathtools}
\usepackage{booktabs,array,enumitem,microtype,xcolor}
\usepackage{tikz-cd}
\usepackage[colorlinks=true,linkcolor=blue!50!black,citecolor=blue!50!black,urlcolor=blue!50!black]{hyperref}
\newtheorem{theorem}{Theorem}[section]
\newtheorem{proposition}[theorem]{Proposition}
\newtheorem{lemma}[theorem]{Lemma}
\newtheorem{corollary}[theorem]{Corollary}
\newtheorem{conjecture}[theorem]{Conjecture}
\newtheorem{question}[theorem]{Question}
\theoremstyle{definition}
\newtheorem{definition}[theorem]{Definition}
\theoremstyle{remark}
\newtheorem{remark}[theorem]{Remark}
\numberwithin{equation}{section}
\newcommand{\RR}{\mathbb R}
\newcommand{\ZZ}{\mathbb Z}

\newcommand{\volsimp}[1]{\lVert #1\rVert}
\DeclareMathOperator{\lk}{lk}
\DeclareMathOperator{\st}{st}
\DeclareMathOperator{\im}{im}

\DeclareMathOperator{\supp}{supp}
\DeclareMathOperator{\OO}{O}
\DeclareMathOperator{\id}{id}

\newcommand{\PP}{\mathcal P}
\newcommand{\CC}{\mathcal C}
\newcommand{\cstar}[1]{\overline{\operatorname{St}}^{\square}_{\Sigma}(#1)}
\AtBeginDocument{%
  \begingroup
    \normalfont\normalsize
    \global\normalparindent=4\fontdimen2
  \endgroup
  \parindent=\normalparindent
  \listparindent=\normalparindent
}
\setlist{listparindent=\normalparindent}
\setlist[enumerate]{label=\textup{(\roman*)},leftmargin=2.1em}
\setlist[itemize]{leftmargin=1.8em}
\hypersetup{pdftitle={Positivity of simplicial volume and uniqueness for the twelve-vertex flag three-spheres},pdfauthor={Xiaojun Lai and Xueyuan Wan}}

\title[Positivity and uniqueness for twelve-vertex flag $3$-spheres]{Positivity of simplicial volume and uniqueness for the twelve-vertex flag three-spheres}
\author{Xiaojun Lai}
\address{Xiaojun Lai: School of Mathematical Sciences, Chongqing University of Technology, Chongqing 400054, China}
\email{xiaojun\_lai@stu.cqut.edu.cn}

\author{Xueyuan Wan}
\address{Xueyuan Wan: Mathematical Sciences Research Center, Chongqing University of Technology,
Chongqing 400054, China}
\email{xwan@cqut.edu.cn}

\date{}
\makeatletter
\@namedef{subjclassname@2020}{\textup{2020} Mathematics Subject Classification}
\makeatother

\begin{document}
\begin{abstract}
We study the twelve-vertex flag triangulation $T_{12}$ of $S^3$ and the
Davis $4$-manifold $M(T_{12})$ associated with it. We prove that
$T_{12}$ is, up to isomorphism, the unique flag triangulation of $S^3$
with twelve vertices in which every edge is contained in a square; this proves Milizia's Conjecture~6.31. We also prove a general criterion for
the tangent Euler class of a Davis $4$-manifold to be bounded, and we
verify it for $T_{12}$ by means of an explicit symmetric form. Since the
Euler class takes the value $256$ on the fundamental class of
$M(T_{12})$, the duality between bounded cohomology and simplicial
volume gives $\|M(T_{12})\|>0$. This answers Milizia's
Question~1.7.
\end{abstract}
\maketitle

\section{Introduction}\label{sec:introduction}

For an oriented closed $n$-manifold $X$, the simplicial volume
$\|X\|$ is the infimum of the $\ell^1$-norms of the real singular
cycles representing the fundamental class. It is a homotopy invariant,
although it is defined by means of chains and not by homology classes
alone. The dual notion is bounded cohomology: if a class in
$H^n(X;\mathbb R)$ is represented by a bounded cocycle and takes a
nonzero value on the fundamental class, then $\|X\|>0$ by the duality
inequality~\cite[Section~1.1, pp.~16--17]{Gromov1982}. To apply this
inequality to a nonzero characteristic number, one needs a bounded
cocycle representing the corresponding characteristic class.

Flag triangulations of spheres provide closed aspherical manifolds to
which this method can be applied. To a flag triangulation $K$ of
$S^{3}$, the reflection-group construction of Davis associates a closed
aspherical $4$-manifold
$M(K)$~\cite[Section~1.2 and Propositions~1.2.2--1.2.3]{Davis2008}.
If $K$ has $m$ vertices,
then
\begin{equation}\label{eq:intro-euler}
\chi(M(K))=2^{m-4}\gamma_2(K),
\qquad \gamma_2(K)=f_1(K)-5f_0(K)+16,
\end{equation}
and $\gamma_2(K)$ is nonnegative for every flag $3$-sphere $K$ by a
theorem of Davis and Okun~\cite[Theorem~11.2.1]{DavisOkun2001}.
When $\gamma_2(K)>0$, formula~\eqref{eq:intro-euler} gives
$\chi(M(K))>0$. To deduce positive simplicial volume from the duality
inequality, it suffices to show that the tangent Euler class of $M(K)$
is bounded.

In this setting,
Milizia~\cite[Definition~1.2, Proposition~4.2, and Corollary~4.4]{Milizia2024}
studies simplicial maps between spheres and the domination relations
which they induce between Davis manifolds. Two small flag $3$-spheres,
denoted by $T_{10}$ and $T_{12}$, appear there as basic examples. The
manifold $M(T_{10})$ is homeomorphic to the product of two closed
surfaces of genus five, and hence has positive simplicial
volume~\cite[Section~1.2, p.~5]{Milizia2024}. The corresponding question
for $T_{12}$ is less obvious. The complex $T_{12}$ is the twelve-vertex
flag $3$-sphere $\Delta_{12,33}$ of
Venturello~\cite[Theorem~3.5]{Venturello2022}; for this identification
see~\cite[Section~6.3, footnote~3]{Milizia2024}.

Following Milizia~\cite[Definition~2.6]{Milizia2024}, we call a full
subcomplex of a flag complex a \emph{square} if it is isomorphic to the
cycle $C_4$. Milizia proposed the following
conjecture~\cite[Conjecture~6.31]{Milizia2024}.

\begin{conjecture}
$T_{12}$ is, up to isomorphism, the only flag triangulation of $S^3$
with twelve vertices in which every edge is contained in a square.
\end{conjecture}

Our first main result confirms this conjecture.

\begin{theorem}\label{thm:uniqueness}
$T_{12}$ is, up to isomorphism, the only flag triangulation of $S^3$
with twelve vertices in which every edge is contained in a square.
\end{theorem}

We begin with a vertex $v$ of maximum degree. We show that exactly
three other vertices are not adjacent to $v$, and that these three
vertices span a triangle. Each edge of this triangle determines an
induced path in the link $\operatorname{lk}(v)$, which is a flag
$2$-sphere. The three paths have the same two endpoints and disjoint
interiors. Using the fact that edges cannot cross in this sphere, we
restrict the possible arrangements of the remaining vertices and edges.
Finally, the links of the two common endpoints, together with the
bounds on vertex degrees, determine all remaining edges. We use the
Davis--Okun inequality~\cite[Theorem~11.2.1]{DavisOkun2001} and
Chudnovsky and Nevo's result that every flag $2$-sphere other than the
octahedral sphere has an edge contained in no
square~\cite[Lemma~3.6]{ChudnovskyNevo2020}.

Concerning the associated Davis manifold, Milizia asked whether
$M(T_{12})$ has positive simplicial volume~\cite[Question~1.7]{Milizia2024}.

\begin{question}\label{question}
Does $M(T_{12})$ have positive simplicial volume?
\end{question}

We answer this question by realizing the tangent Euler class as a
primary characteristic class of a flat principal bundle with connected
semisimple structure group. The construction uses the cubical structure
of $M(K)$ and a symmetric form associated with $K$, and gives the
following criterion for boundedness. A nondegenerate symmetric form is
said to have \emph{inertia} $(p,q)$ if its positive and negative indices
are $p$ and $q$, respectively; this is the usual
notion~\cite[Definition~4.5.6]{HornJohnson2013}, with the zero index
omitted.

\begin{theorem}\label{thm:criterion}
Let $K$ be a flag $3$-sphere with $m$ vertices. Suppose that there is a
nondegenerate symmetric form $B$ on $\RR^m$ such that, in the basis
indexed by $V(K)$,
\begin{equation}\label{eq:criterion-assumptions}
B_{ii}=1,\qquad
B_{ij}=0\quad\text{for each edge }ij\text{ of }K,\qquad
B_{ij}\le0\quad\text{if }i\ne j,
\end{equation}
and that $B$ has inertia $(4,m-4)$. Then the real Euler class
$e(\tau M(K))$ is bounded. Consequently, if $\chi(M(K))\ne0$, then
$\volsimp{M(K)}>0$.
\end{theorem}

The proof of Theorem~\ref{thm:criterion} first relates the form $B$
to the local topology of the Davis manifold. The conditions on $B$
imply that orthogonal projection maps the coordinate cone over $K$
homeomorphically onto any positive four-dimensional subspace.
Using these projections on the Davis cubical complex, we identify
the tangent microbundle of $M(K)$ with the microbundle underlying a
positive rank-four subbundle of a flat rank-$m$ vector bundle.
Milnor's tangent microbundle~\cite[Lemma~2.1]{Milnor1964} describes
the local tangent structure of a topological manifold; Kister's
theorem~\cite[Theorem~2]{Kister1964} represents it by a topological
$\mathbb R^4$-bundle.

We then show that the reflection representation of the commutator
subgroup takes values in $\operatorname{SO}_0(B)$, the identity
component of $\operatorname{SO}(B)$. This connected semisimple Lie
group preserves the orientations of positive four-planes. The Euler
class of the positive subbundle is therefore a primary characteristic
class of the resulting flat principal $\operatorname{SO}_0(B)$-bundle.
The boundedness theorem of Chatterji, de Cornulier, Mislin, and
Pittet~\cite[Theorem~1.1]{ChatterjiCornulierMislinPittet2013} makes
this Euler class bounded. The positive subbundle need not itself be
flat; the theorem applies to the flat principal bundle from which its
Euler class is obtained.

For $T_{12}$ we take the symmetric form whose diagonal entries are $1$,
whose entries on edges are $0$, and whose entries on nonedges are $-4$.
An exact computation of principal minors shows that its inertia is
$(4,8)$. Since $\gamma_2(T_{12})=1$, we have $\chi(M(T_{12}))=256$
by~\eqref{eq:intro-euler}. Therefore Theorem~\ref{thm:criterion} gives
an affirmative answer to Question~\ref{question}.

\begin{corollary}\label{thm:positive}
We orient $M(T_{12})$ by transporting the orientation constructed in
the proof of Theorem~\ref{thm:criterion} through the microbundle
isomorphism in Proposition~\ref{prop:tangent-identification}. With this
orientation, the Euler class of the tangent microbundle of $M(T_{12})$
belongs to
the image of the comparison map
\[
H_b^4(M(T_{12});\RR)\longrightarrow H^4(M(T_{12});\RR).
\]
Its value on the fundamental class is $256$. In particular,
\[
\volsimp{M(T_{12})}>0.
\]
\end{corollary}

The proofs of Theorem~\ref{thm:uniqueness} and
Corollary~\ref{thm:positive} are independent of each other. The proof of
the theorem is based on the combinatorics of links in flag complexes and
on planar separation; that of the corollary uses the explicit form, the
Davis construction, the theory of microbundles, and bounded
characteristic classes. In Section~\ref{sec:preliminaries} we fix the
combinatorial, topological, and cohomological conventions.
Section~\ref{uni:section} is devoted to the proof of the classification
theorem. In Section~\ref{sec:euler} we prove the reflection criterion,
and in Section~\ref{sec:application} we verify its hypotheses for the
form on $T_{12}$, which completes the proof of the corollary.

\section{Preliminaries}\label{sec:preliminaries}

In this section we collect the definitions and known results used in the
proofs. We first treat graphs, flag complexes, and links, and then
introduce $\gamma_2$ and the combinatorial model of $T_{12}$. Next we
describe Davis manifolds and their right-angled Coxeter groups, together
with the relation between $\gamma_2$ and the Euler characteristic.
Finally, we recall the estimate of simplicial volume by bounded
cohomology, and the facts on tangent microbundles, Euler classes, and
flat bundles which are needed for the proof of boundedness.

\subsection{Graphs, flag complexes, and links}
We fix the notation for the combinatorial part of the paper, and prove
some elementary facts on adjacency, fullness, and links for later
reference. For the terminology on simplicial complexes we refer
to~\cite[Sections~A.2--A.3 and~A.6]{Davis2008}.
All graphs and simplicial complexes are assumed to be finite unless
otherwise stated.

\begin{definition}\label{def:graph-paths}
A \emph{simple graph} $G=(V,E)$ consists of a vertex set $V$ and an
edge set $E\subseteq\binom{V}{2}$; thus it has neither loops nor
multiple edges. We write $uv\in E(G)$ for $\{u,v\}\in E$ and put
\[
N_G(u)=\{v\in V:uv\in E(G)\},\qquad \deg_G(u)=|N_G(u)|.
\]
For $A\subseteq V$, the \emph{induced subgraph} $G[A]$ has the vertex
set $A$ and contains exactly those edges of $G$ whose endpoints lie in
$A$. A vertex of degree one is called a \emph{leaf}. The \emph{degree
sequence} is the list of all vertex degrees, with repetitions; the
notation $a^r b^s$ means that degree $a$ occurs $r$ times and degree $b$
occurs $s$ times.

A \emph{simple path} $v_0-v_1-\cdots-v_k$ has distinct vertices and the
edges $v_{i-1}v_i$ for $1\le i\le k$. Its \emph{internal vertices}
are $v_1,\ldots,v_{k-1}$. It is said to be \emph{induced} if no two
nonconsecutive vertices are adjacent. A \emph{cycle} $C_k$, $k\ge3$,
has the vertices $v_0,\ldots,v_{k-1}$ and only the edges
$v_iv_{i+1}$, the indices being taken modulo $k$. Paths with the same
endpoints are said to be \emph{internally vertex-disjoint} if their
sets of internal vertices are pairwise disjoint.
\end{definition}

\begin{lemma}\label{lem:leaf-connected}
Let $a$ be a leaf of $G$, with $N_G(a)=\{b\}$. If $A\subseteq V(G)$,
$a\in A$, $b\notin A$, and $G[A]$ is connected, then $A=\{a\}$.
\end{lemma}
\begin{proof}
Since the only neighbor of $a$ in $G$ does not belong to $A$, the
vertex $a$ has no neighbor in $G[A]$. A connected graph containing an
isolated vertex has no other vertices.
\end{proof}

\begin{definition}
A \emph{simplicial complex} $K$ on $I$ is a collection of subsets of
$I$ which is closed under taking subsets and contains all singletons.
Its elements are called \emph{simplices}; a simplex with $r+1$ vertices
has dimension $r$. The empty simplex is included. We write $V(K)=I$,
denote by $K^{(1)}$ the graph of vertices and edges, and denote by
$f_i(K)$ the number of $i$-dimensional simplices; thus $f_{-1}(K)=1$.
For $A\subseteq I$, the \emph{full}, or \emph{induced}, subcomplex is
\[
K[A]=\{\sigma\in K:\sigma\subseteq A\}.
\]
Thus $K[A]^{(1)}=K^{(1)}[A]$. We use $N_K(u)$ and $\deg_K(u)$ for
the corresponding notions in the graph $K^{(1)}$.
A \emph{simplicial isomorphism} $\varphi:K\to L$ is a bijection
$\varphi:V(K)\to V(L)$ such that, for every $A\subseteq V(K)$,
\[
A\in K\quad\Longleftrightarrow\quad\varphi(A)\in L.
\]
We write $K\cong L$ if such a bijection exists. Accordingly, ``up to
isomorphism'' refers to the simplices, and not merely to the
homeomorphism type of the realizations; this is the convention of
\cite[Introduction, p.~3]{Milizia2024}.
\end{definition}

We use the standard geometric realization in coordinates indexed
by $I$:
\begin{equation}\label{eq:geometric-realization}
|K|=\left\{x\in\RR^I:x_i\ge0,\ \sum_{i\in I}x_i=1,
\ \supp(x)\in K\right\},\quad
\supp(x)=\{i\in I:x_i>0\}.
\end{equation}
The numbers $x_i$ are called the \emph{barycentric coordinates} of $x$.
In particular, $x\in|K[A]|$ if and only if $\supp(x)\subseteq A$.
Each simplex is realized as the convex hull of the coordinate vectors
corresponding to its vertices, and two distinct simplices intersect in
their common face. The relative interior of an edge, i.e., the edge
with its two endpoints removed, contains no vertex of the complex. A
complex is connected if and only if its $1$-skeleton is connected,
since each simplex is connected and any two vertices of a simplex are
joined by an edge. For these conventions see
\cite[Example~A.1.4 and Section~A.2, equation~(A.6)]{Davis2008}.

\begin{definition}\label{def:clique-complex}
A \emph{clique} in a simple graph $G$ is a set of vertices any two
distinct members of which are adjacent; the empty set and singletons
are regarded as cliques. The \emph{clique complex} of $G$ is
\[
\operatorname{Cliq}(G)=\{A\subseteq V(G): A\text{ is a clique in }G\}.
\]
A simplicial complex $K$ is said to be \emph{flag} if every clique in
$K^{(1)}$ is a simplex of $K$; see~\cite[Definition~A.3.5 and
Example~A.3.8]{Davis2008}.
\end{definition}

\begin{proposition}[{\cite[Lemma~A.3.7 and Example~A.3.8]{Davis2008}}]\label{prop:flag-clique}
For every graph $G$, $\operatorname{Cliq}(G)$ is flag and
$\operatorname{Cliq}(G)^{(1)}=G$. For every simplicial complex $K$,
\[
K\text{ is flag}\quad\Longleftrightarrow\quad K=\operatorname{Cliq}(K^{(1)}).
\]
Every full subcomplex of a flag complex is flag. A graph isomorphism
between the $1$-skeleta of two flag complexes induces a simplicial
isomorphism between them.
\end{proposition}

\begin{definition}
For complexes $K,L$ on disjoint vertex sets, their \emph{join} is
\[
K*L=\{\sigma\cup\tau:\sigma\in K,\ \tau\in L\}.
\]
The \emph{suspension} of $K$ is $K*S^0$, where $S^0$ consists of two
vertices and no edge. The \emph{link} and the \emph{closed star} of
$\sigma\in K$ are
\[
\lk_K(\sigma)=\{\tau\in K:\tau\cap\sigma=\varnothing,
\ \tau\cup\sigma\in K\},\qquad
\st_K^{\mathrm{cl}}(\sigma)=\sigma*\lk_K(\sigma),
\]
where, in the join, $\sigma$ denotes the simplex together with all its
faces.
\end{definition}

\begin{lemma}\label{lem:link-identities}
Let $K$ be a simplicial complex and let $\sigma\in K$.
\begin{enumerate}
\item If $\sigma\cap\tau=\varnothing$ and $\sigma\cup\tau\in K$, then
\begin{equation}\label{eq:iterated-link}
\lk_K(\sigma\cup\tau)=\lk_{\lk_K(\sigma)}(\tau).
\end{equation}
In particular, for an edge $ab$,
$\lk_K(ab)=\lk_{\lk_K(a)}(b)$.
\item If $K$ is flag and $\sigma\ne\varnothing$, put
\[
N_K(\sigma)=\bigcap_{u\in\sigma}N_K(u).
\]
Then $N_K(\sigma)\cap\sigma=\varnothing$ and
\begin{equation}\label{eq:link-full-neighbors}
\lk_K(\sigma)=K[N_K(\sigma)].
\end{equation}
Hence links in a flag complex are full and flag.
\item If $K$ is flag, then
$\st_K^{\mathrm{cl}}(u)=K[\{u\}\cup N_K(u)]$ for every vertex $u$.
\end{enumerate}
\end{lemma}
\begin{proof}
(i) A simplex $\eta$ belongs to the right-hand side if and only if
$\eta\cap\tau=\varnothing$, $(\eta\cup\tau)\cap\sigma=\varnothing$,
and $\eta\cup\tau\cup\sigma\in K$. These are exactly the conditions
for $\eta$ to belong to the left-hand side.
Compare~\cite[Lemma~A.6.3]{Davis2008}.

(ii) Every vertex of the link is adjacent to every vertex of $\sigma$.
Conversely, if a vertex is adjacent to all vertices of $\sigma$, then
its union with $\sigma$ belongs to $K$, since $K$ is flag. Let
$\eta\in K[N_K(\sigma)]$. Then any two vertices of $\eta\cup\sigma$
are joined by an edge: the pairs inside $\eta$ and inside $\sigma$
because these are simplices, and the remaining pairs by the definition
of $N_K(\sigma)$. Since $K$ is flag, $\eta\cup\sigma\in K$, so that
$\eta$ lies in the link. The reverse inclusion follows from the
definition. The link is flag by Proposition~\ref{prop:flag-clique}.

(iii) A simplex on $\{u\}\cup N_K(u)$ either contains $u$, or all its
vertices are adjacent to $u$; in the latter case its union with $u$ is
a simplex, since $K$ is flag. In either case it belongs to the closed
star. The other inclusion is obvious from the definition of the star.
\end{proof}

\begin{definition}\label{def:flag-square}
A \emph{flag $d$-sphere} is a flag complex whose realization is
homeomorphic to $S^d$. Following~\cite[Definition~2.6]{Milizia2024},
a \emph{square} in $K$ is a full subcomplex isomorphic to $C_4$.
Thus a square $a-b-c-d-a$ has four distinct vertices and the four
edges $ab,bc,cd,da$, and contains neither of the diagonals $ac,bd$.
Equivalently, these vertices span an induced $4$-cycle in $K^{(1)}$.
\end{definition}

The next lemma is a low-dimensional consequence of a criterion for
homology manifolds.
\begin{lemma}[{\cite[Lemma~10.4.6]{Davis2008}}]\label{lem:sphere-links}
Let $K$ triangulate $S^d$, where $d=2$ or $3$. If $\sigma$ is an
$r$-simplex with $0\le r<d$, then $\lk_K(\sigma)$ triangulates
$S^{d-r-1}$. In particular:
\begin{enumerate}
\item in a triangulated $2$-sphere, every vertex link is a cycle and
 every edge belongs to exactly two triangles;
\item in a triangulated $3$-sphere, vertex links are $2$-spheres,
 edge links are cycles, and each triangle belongs to exactly two
 tetrahedra, whose additional vertices form its link $S^0$.
\end{enumerate}
If $K$ is flag, all these links are full and flag, and each cycle
appearing as a link has at least four vertices.
\end{lemma}
\begin{proof}
By the cited criterion, every link is a generalized homology sphere,
i.e., a homology manifold with the integral homology of a sphere; here
its dimension is at most two. A $0$-dimensional homology sphere consists
of two points. A $1$-dimensional generalized homology sphere is a
connected graph in which every vertex has degree two, and hence it is a
cycle. Consequently, in dimension two the vertex links are circles, so
that the link is a closed connected surface with the integral homology
of $S^2$; by the classification of closed surfaces it is $S^2$.
The assertions on incidences follow by considering the links of
simplices of codimension one and two. If $K$ is flag, the links are full
and flag by Lemma~\ref{lem:link-identities}, and a flag cycle cannot
have three vertices.
\end{proof}

\begin{lemma}\label{lem:flag-links}
A flag $2$-sphere $K$ has
\[
\deg_K(u)\ge4\quad(u\in V(K)),\qquad
f_1(K)=3f_0(K)-6,\qquad f_2(K)=2f_0(K)-4.
\]
In particular $f_0(K)\ge6$, with equality if and only if
$K\cong(S^0)^{*3}$, the octahedral sphere.
\end{lemma}
\begin{proof}
The link of $u$ is a cycle with $\deg_K(u)$ vertices, and by
Lemma~\ref{lem:sphere-links} it has at least four vertices.
Counting incidences between edges and triangles, we get $3f_2=2f_1$.
Combining this with $f_0-f_1+f_2=2$, we obtain the displayed face
numbers. Since the sum of the vertex degrees counts each edge twice, we
have
\[
6f_0-12=2f_1=\sum_{u\in V(K)}\deg_K(u)\ge4f_0,
\]
and hence $f_0\ge6$.

Suppose that $f_0=6$. Then the degree sum is $24$, and since every
vertex has degree at least four, all six vertices have degree exactly
four. Thus each vertex is joined to four of the other five vertices and
is not joined to the remaining one. If $a$ and $b$ are not joined, then
each of them is the unique nonneighbor of the other. Hence the vertices
form three disjoint pairs, such that there is no edge within a pair and
any two vertices from different pairs are joined.

Therefore a set of vertices is a clique if and only if it contains at
most one vertex from each pair. Since $K$ is flag, these sets are
exactly the simplices of $K$ by Proposition~\ref{prop:flag-clique}.
Each pair, with no edge between its two vertices, is a copy of $S^0$,
and the simplices of the join of the three copies are precisely the
sets containing at most one vertex from each pair. Hence the join is
$K$, and $K\cong(S^0)^{*3}$, the boundary of an octahedron. Conversely,
this boundary is a flag $2$-sphere with six vertices.
\end{proof}

\begin{theorem}[{\cite[Lemma~3.6]{ChudnovskyNevo2020}}]
\label{thm:octahedral}
Every flag $2$-sphere other than the octahedral sphere has an edge
contained in no square.
\end{theorem}

\subsection{Face enumeration and the quantity \texorpdfstring{$\gamma_2$}{gamma2}}
The $h$-polynomial expresses the face numbers in the form used in the
enumerative theory of simplicial complexes; see
Stanley~\cite[Section~II.2, p.~58]{Stanley1996}. The $\gamma$-expansion
below is that of Gal~\cite[Proposition~2.1.1 and Definition~2.1.4]{Gal2005};
for flag spheres it is also used by
Chudnovsky and Nevo~\cite[Section~2.2]{ChudnovskyNevo2020}.

\begin{definition}[{\cite[Section~II.2, p.~58]{Stanley1996}}]
For a $(d-1)$-dimensional simplicial complex $K$, its \emph{$h$-polynomial}
is
\[
h_K(t)=\sum_{i=0}^{d}f_{i-1}(K)t^i(1-t)^{d-i}.
\]
\end{definition}

\begin{theorem}[{\cite[Corollary~II.5.2, p.~65, and Section~II.6, p.~67]{Stanley1996}}]
\label{thm:dehn-sommerville}
If $K$ is a simplicial $(d-1)$-sphere, then
$h_K(t)=t^d h_K(t^{-1})$. Equivalently, the coefficients of
$h_K(t)=\sum_{i=0}^{d}h_i t^i$ satisfy $h_i=h_{d-i}$.
\end{theorem}

The symmetry $h_i=h_{d-i}$ for simplicial spheres is also stated in
\cite[Section~1.2]{Gal2005}.

\begin{definition}[{\cite[Proposition~2.1.1 and Definition~2.1.4]{Gal2005}}]
For a simplicial $(d-1)$-sphere $K$, the coefficients $\gamma_i(K)$ are
defined by the expansion
\[
h_K(t)=\sum_{i=0}^{\lfloor d/2\rfloor}
\gamma_i(K)t^i(1+t)^{d-2i}.
\]
This expansion exists and is unique, because the polynomials
$t^i(1+t)^{d-2i}$ form a basis of the polynomials with the symmetry of
Theorem~\ref{thm:dehn-sommerville}; moreover, $\gamma_0(K)=1$.
\end{definition}

\begin{lemma}[{\cite[Lemma~2.2]{LabbeNevo2017}}]\label{lem:gamma-formula}
For a simplicial $3$-sphere,
\[
h_K(t)=(1+t)^4+\gamma_1(K)t(1+t)^2+\gamma_2(K)t^2,
\]
where
\begin{equation}\label{eq:gamma-definition}
\gamma_1(K)=f_0(K)-8,\qquad
\gamma_2(K)=f_1(K)-5f_0(K)+16.
\end{equation}
\end{lemma}
In the cited lemma, $f_1$ and $f_2$ denote the numbers of vertices and
edges, respectively; with $d=4$ and our indexing one obtains the
formulas above.

Following Milizia~\cite[Section~6, p.~25]{Milizia2024}, we call
$\gamma_2$ a \emph{quantity}. It is invariant under simplicial
isomorphisms, but it is not an invariant of the underlying topological
space: $(S^0)^{*4}$ and $T_{12}$ both triangulate $S^3$, whereas their
values of $\gamma_2$ are $0$ and $1$, respectively. The latter value is
computed in Lemma~\ref{lem:T12-counts}.

\begin{theorem}[{\cite[Theorem~11.2.1, p.~52]{DavisOkun2001}}]
\label{thm:davis-okun}
If $K$ is a flag $3$-sphere, then $\gamma_2(K)\ge0$.
\end{theorem}
In the notation of~\cite[p.~8]{DavisOkun2001}, where
$\kappa(K)=\sum_{i=-1}^{3}(-1/2)^{i+1}f_i(K)$, the definitions above
give
\[
16\kappa(K)=h_K(-1)=\gamma_2(K).
\]

\subsection{The sphere $T_{12}$}
\begin{definition}[{\cite[Equation~(4) and Theorem~3.5]{Venturello2022}; \cite[Figure~12 and footnote~3]{Milizia2024}}]
We denote by $T_{12}$ the clique complex of the graph in
Table~\ref{tab:T12}, with vertex set $\{0,\ldots,11\}$ labeled as
in~\cite[Figure~12]{Milizia2024}. Thus the simplices of $T_{12}$ are
exactly the cliques of this graph. It is Venturello's flag $3$-sphere
$\Delta_{12,33}$, which is identified with $T_{12}$ in
\cite[Section~6.3, footnote~3]{Milizia2024}.
\end{definition}
Every edge of $T_{12}$ is contained in a square; see
\cite[proof of Theorem~3.6]{Venturello2022}.

\begin{table}[htbp]
\centering
\caption{The graph of $T_{12}$, labeled as in
\cite[Figure~12]{Milizia2024}.}
\label{tab:T12}
\renewcommand{\arraystretch}{1.08}
\begin{tabular}{@{}c l@{}}\toprule
$i$ & $N_{T_{12}}(i)$\\\midrule
0 & $1,2,3,4,6,8,9,11$\\
1 & $0,2,4,6,7,8,9,10$\\
2 & $0,1,3,4,5,8,10,11$\\
3 & $0,2,4,5,6,7,9,11$\\
4 & $0,1,2,3,5,7,9,10$\\
5 & $2,3,4,6,7,8,10,11$\\
6 & $0,1,3,5,7,8,9,11$\\
7 & $1,3,4,5,6,8,9,10$\\
8 & $0,1,2,5,6,7,10,11$\\
9 & $0,1,3,4,6,7$\\
10 & $1,2,4,5,7,8$\\
11 & $0,2,3,5,6,8$\\\bottomrule
\end{tabular}
\end{table}

\begin{lemma}\label{lem:T12-counts}
The sphere $T_{12}$ has degree sequence $6^3 8^9$: vertices
$9,10,11$ have degree six, and vertices $0,\ldots,8$ have degree eight.
Its face numbers and $\gamma_2$ are
\begin{equation}\label{eq:T12-fvector}
(f_0,f_1,f_2,f_3)=(12,45,66,33),\qquad \gamma_2(T_{12})=1.
\end{equation}
\end{lemma}
\begin{proof}
From the degrees in the table we obtain
\[
2f_1=\sum_{i=0}^{11}\deg_{T_{12}}(i)=3\cdot6+9\cdot8=90.
\]
Since every tetrahedron has four triangles and every triangle lies in
two tetrahedra, we have $4f_3=2f_2$. Euler's formula for $S^3$ gives
$f_0-f_1+f_2-f_3=0$. Hence
\[
f_3=f_1-f_0=33,\qquad f_2=66,\qquad
\gamma_2(T_{12})=45-5\cdot12+16=1.
\]
\end{proof}

\subsection{Davis manifolds and right-angled Coxeter groups}
\begin{definition}[{\cite[Definitions~A.2.12, A.6.1, and~I.3.1]{Davis2008}}]
In this paper a \emph{cubical complex} is a space obtained by gluing
Euclidean cubes along faces by isometries, in such a way that the cubes
are embedded and any two of them intersect in a face or not at all. The
\emph{link of a vertex} has one vertex for each incident edge, and one
simplex for each collection of incident edges contained in a common
cube. Thus $r$ such edges in an $r$-cube determine an $(r-1)$-simplex
of the link.
\end{definition}

\begin{definition}[{\cite[Definition~3.1]{Milizia2024}}]
Let $K$ be a finite flag complex on $I$, with $|I|=m$.
Define the cubical subcomplex
\begin{equation}\label{eq:davis-cubical-model}
M(K)=\bigcup_{\sigma\in K}[-1,1]^\sigma\times\{-1,1\}^{I\setminus\sigma}
\ \subseteq[-1,1]^I.
\end{equation}
In the term indexed by $\sigma$, the coordinates in $\sigma$ vary in
$[-1,1]$, while each remaining coordinate is fixed at $+1$ or $-1$.
More explicitly, a simplex $\sigma\in K$ and a choice of signs
$\epsilon\in\{-1,1\}^{I\setminus\sigma}$ determine the cube
\[
Q(\sigma,\epsilon)=\{x\in[-1,1]^I:x_i=\epsilon_i\ (i\notin\sigma)\}.
\]
Thus $M(K)$ consists of the coordinate faces whose sets of free
coordinates are simplices of $K$. For each $\sigma\in K$, the union
in~\eqref{eq:davis-cubical-model} contains $2^{m-|\sigma|}$ such
$|\sigma|$-dimensional faces. If $K$ is a flag $3$-sphere, we call
$M(K)$ the \emph{Davis manifold} of $K$.
\end{definition}

\begin{definition}[{\cite[Section~1.2, pp.~11--12]{Davis2008}}]
For a finite flag complex $K$ on $I$, set
\[
W_K=\langle s_i\ (i\in I)\mid s_i^2=1,
\ s_is_j=s_js_i\text{ for each edge }ij\text{ of }K\rangle.
\]
For $\sigma\in K$, put
$W_\sigma=\langle s_i:i\in\sigma\rangle\cong(\ZZ/2)^{|\sigma|}$.
The \emph{Davis cubical complex} $\Sigma_K$ has the elements of $W_K$
as its vertices, and one $|\sigma|$-cube with vertex set $wW_\sigma$
for each coset $wW_\sigma$; the face relation is given by inclusion of
cosets. The group $W_K$ acts on $\Sigma_K$ by left multiplication; see
also \cite[Proposition~7.3.4 and Example~7.3.6]{Davis2008}.
\end{definition}
The edges at $w$ join $w$ to the vertices $ws_i$. If $ij$ is an edge of
$K$, the four vertices $w,ws_i,ws_is_j,ws_j$ are the vertices of a
$2$-dimensional cube. More generally, the $|\sigma|$ commuting
involutions $s_i$, $i\in\sigma$, give the $2^{|\sigma|}$ vertices of a
cube. For involutions the relation $s_is_j=s_js_i$ is equivalent to
$(s_is_j)^2=1$; it is the relation satisfied by reflections in
perpendicular hyperplanes, which explains the term ``right-angled.'' No
relation is imposed for a nonedge. Thus $W_K$ is defined as an abstract
group; a particular linear reflection representation of it is
constructed in Section~\ref{sec:euler}. The link of every vertex of
$\Sigma_K$ is naturally isomorphic to $K$, and $\dim\Sigma_K=\dim K+1$.
Since each vertex belongs to exactly one cube of each type $\sigma\in K$,
the complex $\Sigma_K$ is locally finite.

\begin{theorem}[{\cite[Section~1.2, pp.~10--13, and Theorem~12.2.1(i)]{Davis2008}; \cite[Lemma~4.3]{Milizia2024}}]
\label{thm:davis-manifold}
For a finite flag complex $K$, the complex $\Sigma_K$ is locally
finite and contractible. Define
\[
\varepsilon:W_K\longrightarrow(\ZZ/2)^I,\qquad s_i\longmapsto e_i,
\]
where $e_i$ is the $i$th standard generator. Then
\begin{equation}\label{eq:commutator}
\Gamma_K=\ker\varepsilon=[W_K,W_K]
\end{equation}
acts freely, properly, and cocompactly on $\Sigma_K$, and
$\Sigma_K/\Gamma_K\cong M(K)$ as cubical complexes.
If $K$ is a flag $3$-sphere, $M(K)$ is a closed connected orientable
aspherical topological $4$-manifold, with fundamental group
$\Gamma_K$ and universal cover $\Sigma_K$.
\end{theorem}
Here the action is \emph{free} if only the identity fixes a point, and
\emph{cocompact} if the quotient is compact; a space is
\emph{aspherical} if its universal cover is contractible. For this
locally finite cubical action, \emph{proper} means that
for every compact subset $C\subseteq\Sigma_K$ the set
$\{\gamma\in\Gamma_K:\gamma C\cap C\ne\varnothing\}$
is finite; for an equivalent formulation in this setting,
see~\cite[Definition~5.1.5]{Davis2008}. This notion concerns the group
action and should not be confused with properness of a metric space.

\begin{proposition}[{\cite[Lemma~3.3 and Section~6, p.~25]{Milizia2024}}]
\label{prop:davis-euler}
If $K$ is a flag $3$-sphere with $m$ vertices, then
\begin{equation}\label{eq:euler-Davis}
\chi(M(K))=\sum_{r=0}^4(-1)^r2^{m-r}f_{r-1}(K)
=2^{m-4}\gamma_2(K).
\end{equation}
\end{proposition}
The first equality in~\eqref{eq:euler-Davis} is obtained by counting
the cubes of $M(K)$, and the second follows from the expansion in
Lemma~\ref{lem:gamma-formula}. This formula will be used in
Section~\ref{sec:application}.

\subsection{Simplicial volume and bounded cohomology}
Unless otherwise stated, all cohomology groups below have real
coefficients.

\begin{definition}[{\cite[Section~0.2, p.~8]{Gromov1982}}]
For a real singular chain written as a combination of distinct singular
simplices, set $\|\sum_j a_j\sigma_j\|_1=\sum_j|a_j|$. The
\emph{simplicial volume} of an oriented closed $n$-manifold $M$ is
\[
\volsimp M=\inf\{\|z\|_1:z\in C_n(M;\RR),\ \partial z=0,
\ [z]=[M]\}.
\]
\end{definition}

\begin{definition}[{\cite[Section~2]{Bucher2007}}]
The bounded singular cochains on $X$ form the subcomplex
\[
C_b^q(X;\RR)=\{c\in C^q(X;\RR):
\|c\|_\infty=\sup_{\sigma:\Delta^q\to X}|c(\sigma)|<\infty\}.
\]
Its cohomology is denoted by $H_b^*(X;\RR)$. The inclusion induces the
\emph{comparison map} $H_b^*(X;\RR)\to H^*(X;\RR)$, and an ordinary
cohomology class in its image is said to be \emph{bounded}.
\end{definition}

\begin{proposition}[{\cite[Section~1.1, pp.~16--17]{Gromov1982}}]
\label{prop:duality}
If a bounded singular cocycle $c$ represents
$\alpha\in H^n(M;\RR)$ for an oriented closed $n$-manifold $M$, then
\begin{equation}\label{eq:duality}
|\langle\alpha,[M]\rangle|\le\|c\|_\infty\volsimp M.
\end{equation}
\end{proposition}
This duality estimate will be used in the proof of
Corollary~\ref{thm:positive}.

\subsection{Microbundles and Euler classes}\label{subsec:microbundles}
The tangent bundle used below is the tangent microbundle, which is
defined for topological manifolds. We recall its relation to the smooth
tangent bundle and the definition of the Euler class needed in the
proof.

\begin{definition}[{\cite[Section~2, pp.~54--56]{Milnor1964}}]
A \emph{microbundle} of rank $n$ over $X$ is a diagram
\[
X\xrightarrow{i}E\xrightarrow{p}X,\qquad p\circ i=\id_X,
\]
such that each $x\in X$ has an open neighborhood $U$ for which a
neighborhood of $i(U)$ in $p^{-1}(U)$ is homeomorphic over $U$ to
$U\times\RR^n$ by a homeomorphism carrying $i$ to the zero section.
Two microbundles over $X$ are isomorphic if neighborhoods of their zero
sections are homeomorphic by a map commuting with the projections and
the sections. Thus the isomorphism class depends only on the germ near
the section. Every vector bundle determines in this way an underlying
microbundle.
\end{definition}

\begin{definition}[{\cite[Lemma~2.1]{Milnor1964}}]
For a topological $n$-manifold $M$, its \emph{tangent microbundle} is
\begin{equation}\label{eq:tangent-microbundle}
\tau M=\bigl(M\xrightarrow{\Delta}M\times M
\xrightarrow{\operatorname{pr}_1}M\bigr),\qquad \Delta(x)=(x,x).
\end{equation}
\end{definition}
To see that $\tau M$ is a microbundle, note that a coordinate chart
$\phi:U\to\RR^n$ gives the local map
$(x,y)\mapsto(x,\phi(y)-\phi(x))$, whose second component, the
displacement from $x$, vanishes on the diagonal. Under a coordinate
change $h$ the second component becomes
$h(\phi(x)+v)-h(\phi(x))$, which is a local homeomorphism fixing zero.
No differentiability is involved.

\begin{theorem}[{\cite[Theorem~2.2]{Milnor1964}}]
If $M$ is a smooth paracompact manifold, the underlying microbundle
of its tangent vector bundle $TM$ is isomorphic to $\tau M$.
\end{theorem}
If $M$ is given a Riemannian metric, the map
$(x,v)\mapsto(x,\exp_x v)$, restricted to suitable neighborhoods of the
zero section and of the diagonal, gives this isomorphism. Only
neighborhoods are compared here; the total spaces $TM$ and $M\times M$
need not be homeomorphic.

\begin{theorem}[{\cite[Theorem~2]{Kister1964}}]\label{thm:kister}
A rank-$n$ microbundle over a finite-dimensional locally finite
simplicial complex is represented by a topological $\RR^n$-bundle
with zero section, unique up to bundle isomorphism preserving that
section.
\end{theorem}
Every finite cubical complex which occurs as a base space in this paper
has a finite simplicial subdivision. Hence, by
Theorem~\ref{thm:kister}, every microbundle over such a base is
represented by a topological $\RR^n$-bundle. Kister proves the same
statement for topological manifolds~\cite[Corollary~1]{Kister1964}.
The transition maps of a topological $\RR^n$-bundle are homeomorphisms
fixing zero; they need not be linear.

\begin{remark}[{\cite[Theorem~2.1 and Remarks~1--2, pp.~277--279]{Holm1967}}]
Let $p:E\to X$ be a topological $\RR^n$-bundle with zero section $i$,
where $X$ is a compact manifold or a finite CW complex. An orientation
specifies a locally consistent choice of generators of
$H^n(E_x,E_x\setminus\{0\};\ZZ)$. There is a unique integral Thom
class $u\in H^n(E,E\setminus i(X);\ZZ)$ restricting to these generators.
If $j:H^n(E,E\setminus i(X);\ZZ)\to H^n(E;\ZZ)$ is the natural map,
the integral Euler class is $i^*j(u)$; we write $e$ for its image in
$H^n(X;\RR)$. For vector bundles this is the usual construction
\cite[Section~9, pp.~96--98]{MilnorStasheff1974}.
Excision makes it depend only on the germ near the zero section, so
it also defines the Euler class of the oriented microbundle represented
by $E$. Thom and Euler classes are natural under pullback with the
induced orientations.
\end{remark}

\begin{proposition}[{cf. \cite[Section~2, pp.~1401--1402]{Samelson1965}}]
\label{prop:microbundle-euler}
An oriented microbundle isomorphism preserves the Euler class.
For an oriented closed topological $n$-manifold $M$,
\begin{equation}\label{eq:Euler-pairing}
\langle e(\tau M),[M]\rangle=\chi(M).
\end{equation}
\end{proposition}

\subsection{Flat bundles and primary characteristic classes}
We recall the description by classifying maps used in
\cite[Section~3]{Bucher2007}. In this subsection all base spaces
are CW complexes, and all structure groups are Lie groups or discrete
groups.

\begin{definition}[{\cite[Definition~1.4 and Proposition~1.11]{CohenFiberBundles}}]
Let $p:E\to X$ be a bundle and let $f:Z\to X$ be a continuous map.
The \emph{pullback of $E$ along $f$} is
\[
f^*E=\{(z,e)\in Z\times E:f(z)=p(e)\},
\qquad p_f(z,e)=z.
\]
The fiber over $z$ is canonically identified with $E_{f(z)}$. The
defining square is
\[
\begin{tikzcd}
f^*E \arrow[r,"\operatorname{pr}_2"] \arrow[d,"p_f"'] & E\arrow[d,"p"]\\
Z\arrow[r,"f"']&X.
\end{tikzcd}
\]
A \emph{section} of $p$ is a continuous map $s:X\to E$ with
$p\circ s=\id_X$. A section of $f^*E$ has the form
$z\mapsto(z,\widetilde f(z))$, where $p\circ\widetilde f=f$; thus
sections of the pullback correspond to lifts of $f$.
\end{definition}

\begin{definition}[{\cite[Definitions~1.6 and~2.1, and Theorem~2.21]{CohenFiberBundles}}]
Let $G$ be a topological group. A \emph{principal $G$-bundle} is a
locally trivial bundle with fiber $G$ and a right $G$-action such that
the local trivializations are compatible with right multiplication. A
\emph{classifying space} $BG$ is the base of a universal principal
bundle $EG\to BG$; that is, every principal $G$-bundle
$\mathcal H\to X$ is isomorphic to $f_{\mathcal H}^*(EG)$ for a map
$f_{\mathcal H}:X\to BG$, unique up to homotopy. One may choose $EG$
to be contractible; see Milnor's
construction~\cite[Section~3, pp.~432--433]{Milnor1956Universal}.
Thus $[X,BG]$ classifies principal $G$-bundles over $X$ up to
isomorphism.
\end{definition}

\begin{definition}[{\cite[Chapter~1, Section~1.2, pp.~13--14]{CohenFiberBundles}}]
Let $p:\mathcal H\to X$ be a principal $G$-bundle and let $Z$ be a
space with a continuous left $G$-action. The \emph{associated bundle}
is
\[
\mathcal H\times_G Z=(\mathcal H\times Z)/\!\sim,
\qquad (ag,z)\sim(a,gz),
\]
with projection $[a,z]\mapsto p(a)$ and fiber $Z$.
If $Z$ is a vector space on which $G$ acts linearly, this is a vector
bundle. The construction commutes with pullback:
\[
f^*(\mathcal H\times_G Z)\cong(f^*\mathcal H)\times_G Z.
\]
\end{definition}

If $\Lambda$ is a discrete group, $B\Lambda$ is a $K(\Lambda,1)$,
i.e., it has fundamental group $\Lambda$ and contractible universal
cover, and its homotopy type is determined by $\Lambda$
\cite[Example~1B.7 and Theorem~1B.8]{Hatcher2002}. In particular,
since $M(K)$ has fundamental group $\Gamma_K$ and contractible universal
cover $\Sigma_K$, it is itself a model for $B\Gamma_K$.

\begin{definition}[{\cite[Section~3]{Bucher2007}}]
We denote by $G^\delta$ the group $G$ with the discrete topology. A
principal $G$-bundle is said to be \emph{flat} if it admits locally
constant transition functions, that is, continuous transition functions
with values in $G^\delta$. The identity homomorphism $G^\delta\to G$ is
continuous and induces the canonical map $\iota:BG^\delta\to BG$, where
$BG^\delta=B(G^\delta)$; this notation does not mean that $BG$ itself
is given the discrete topology.
\end{definition}
Extending the structure group along $G^\delta\to G$ forgets the
chosen flat structure. Consequently, $\mathcal H$ admits a flat
structure if and only if its classifying map factors up to homotopy as
\[
X\xrightarrow{c_\delta}BG^\delta\xrightarrow{\iota}BG,
\qquad f_{\mathcal H}\simeq\iota\circ c_\delta.
\]

Let $X$ be connected and let $\rho:\pi_1(X)\to G$ be a representation.
Then the quotient
\[
\mathcal G_\rho=(\widetilde X\times G)/\pi_1(X),\qquad
\gamma\cdot(x,g)=(\gamma x,\rho(\gamma)g),
\]
is a flat principal bundle. If $G$ acts linearly on $V$, the associated
vector bundle is $(\widetilde X\times V)/\pi_1(X)$, where $\gamma$ acts
by $(x,v)\mapsto(\gamma x,\rho(\gamma)v)$. Up to conjugacy, the
representation is the \emph{monodromy}, also called the holonomy, of
this flat structure. We do not assume that $\rho(\pi_1(X))$ is discrete
in $G$; since $\pi_1(X)$ is discrete, $\rho$ is continuous both as a
map into $G^\delta$ and as a map into $G$.

For $M=M(K)\simeq B\Gamma_K$, a representation of $\Gamma_K$
therefore gives a map $c_\delta:M\simeq B\Gamma_K\to BG^\delta$, which
is determined up to based homotopy by the
homomorphism~\cite[Proposition~1B.9]{Hatcher2002}. Its composite with
$\iota$ classifies the underlying ordinary bundle.

\begin{definition}[{\cite[Section~3]{Bucher2007}}]
A universal characteristic class $\beta\in H^q(BG;\RR)$ assigns
to $\mathcal G$ the class $f_{\mathcal G}^*\beta$. By a
\emph{primary characteristic class} of flat $G$-bundles we mean a class
in $\im(\iota^*:H^*(BG;\RR)\to H^*(BG^\delta;\RR))$.
For a chosen flat structure, its value on $X$ is
\[
\beta(\mathcal G)=c_\delta^*\iota^*\beta.
\]
\end{definition}

We use the following consequence of the boundedness theorem of
Chatterji, de Cornulier, Mislin, and Pittet.

\begin{theorem}[{\cite[Theorem~1.1]{ChatterjiCornulierMislinPittet2013}}]\label{thm:primary-bounded}
Let $G$ be a connected semisimple Lie group, with its usual topology.
Then
\[
\im\bigl(H^*(BG;\RR)\to H^*(BG^\delta;\RR)\bigr)
\subseteq
\im\bigl(H_b^*(BG^\delta;\RR)\to H^*(BG^\delta;\RR)\bigr).
\]
\end{theorem}
The theorem is stated in terms of the bounded cohomology of the
discrete group $G^\delta$. We use the standard comparison between the
bar-resolution model for this bounded group cohomology and bounded
singular cohomology of the classifying space $BG^\delta$; the latter
model is reviewed in~\cite[Section~2]{Bucher2007}.
The cited theorem characterizes the connected Lie groups for which
this inclusion holds: the derived group of the radical must be simply
connected. The radical is the largest connected normal solvable
subgroup; it is trivial for a semisimple group, so the stated case
follows. Earlier boundedness results for real linear algebraic groups
were proved by Gromov~\cite[Section~1.3, pp.~23--24]{Gromov1982}
and Bucher-Karlsson~\cite[Theorem~7]{Bucher2007}.

If $\zeta$ is a bounded singular cocycle representing
$\iota^*\beta$, then $c_\delta^*\zeta$ represents
$\beta(\mathcal G)$ and
$\|c_\delta^*\zeta\|_\infty\le\|\zeta\|_\infty$.
In Section~\ref{sec:euler} we apply the theorem to the connected
semisimple group $G=\operatorname{SO}_0(B)$. We show that the positive
subbundle $\xi^+$ satisfies
$e(\xi^+)=c_\delta^*\iota^*\mathfrak e^+$ for a class
$\mathfrak e^+\in H^4(BG;\RR)$.
This expresses its Euler class as a primary characteristic class of a
flat principal $G$-bundle, even though $\xi^+$ need not be flat.

\section{The twelve-vertex classification}\label{uni:section}

In this section we prove Theorem~\ref{thm:uniqueness} by determining
the graph of the triangulation. Throughout the section, $T$ is a flag
triangulation of $S^3$ with twelve vertices such that every edge of $T$
belongs to a square. We choose a vertex $v$ of maximum degree and put
\begin{equation}\label{uni:notation}
 d=\deg_T(v),\qquad L=\lk_T(v),\qquad
 A=V(T)\setminus\bigl(\{v\}\cup V(L)\bigr).
\end{equation}
Thus $V(L)=N_T(v)$, $|V(L)|=d$, and $|A|=11-d$.
For each $z\in V(L)$, let
\begin{equation}\label{uni:neighbor-sets}
 S_z=N_T(z)\cap A
\end{equation}
be the set of neighbors of $z$ in $A$. By
Lemmas~\ref{lem:link-identities} and~\ref{lem:sphere-links},
$L$ is a full flag subcomplex of $T$ triangulating $S^2$.
Hence the partition $V(T)=\{v\}\sqcup V(L)\sqcup A$ gives
\begin{equation}\label{uni:degree-decomposition}
 \begin{aligned}
 N_T(z)&=\{v\}\sqcup N_L(z)\sqcup S_z,\\
 \deg_T(z)&=1+\deg_L(z)+|S_z|
 \end{aligned}
 \qquad(z\in V(L)).
\end{equation}

\subsection{Nonneighbors and the maximum degree}

We begin with a connectivity property, which follows from the statement
of Labb\'e and Nevo on the deletion of an equator, applied to the
induced sphere $\lk_K(u)$ of codimension one.

\begin{lemma}[{\cite[Lemma~2.1(iv)]{LabbeNevo2017}}]\label{uni:nonneighbors}
Let $K$ be a flag triangulation of $S^m$, where $m=2$ or $3$, and let $u$
be a vertex. The full subcomplex on the nonneighbors of $u$ is nonempty
and connected.
\end{lemma}
\begin{proof}
Put
\[
\Lambda=\lk_K(u),\qquad
A=V(K)\setminus\bigl(\{u\}\cup N_K(u)\bigr).
\]
Since $K$ is flag, $\Lambda$ is an induced subcomplex, and the deletion of
$\Lambda$ is
\[
K[V(K)\setminus V(\Lambda)]
 =K[\{u\}\cup A]=\{u\}\sqcup K[A].
\]
The vertex $u$ is isolated in this deletion. By
\cite[Lemma~2.1(iv)]{LabbeNevo2017}, the deletion of the induced equator
$\Lambda$ has the reduced homology of $S^0$; in particular, it has exactly
two connected components. Hence $K[A]$ is nonempty and connected.
\end{proof}

The next two lemmas hold for an arbitrary flag $3$-sphere and any
vertex $v$, with $L,A,S_z$ defined as above. Their proofs use neither
the number of vertices, nor the maximality of $\deg_T(v)$, nor the
assumption on squares.

\begin{lemma}\label{uni:labels}
The sets $S_z$ satisfy
\begin{align}
 &S_z\ne\varnothing\quad\text{and}\quad T[S_z]\text{ is connected}
       &&(z\in V(L));\label{uni:label-connected}\\
 &S_x\cap S_y\ne\varnothing
       &&(xy\in L^{(1)});\label{uni:label-edge}\\
 &|S_x\cap S_y\cap S_z|=1
       &&(xyz\in L^{(2)}).\label{uni:label-triangle}
\end{align}
\end{lemma}

\begin{proof}
Fix $z\in V(L)$ and consider the flag $2$-sphere $Q=\lk_T(z)$, which
contains the vertex $v$. Since $Q$ is full, the vertices of $Q$
different from $v$ and nonadjacent to $v$ are precisely the elements
of $N_T(z)\cap A=S_z$. By Lemma~\ref{uni:nonneighbors}, applied to $v$
in $Q$, the complex $Q[S_z]$ is nonempty and connected. Since $Q$ is
full, $Q[S_z]=T[S_z]$, which proves~\eqref{uni:label-connected}.

Let $xyz$ be a triangle of $L$. By Lemma~\ref{lem:sphere-links}, it
belongs to exactly two tetrahedra of $T$. One of them is $vxyz$; we
write the other as $wxyz$. The vertex $w$ is distinct from $v,x,y,z$.
If $w\in V(L)$, then $w,x,y,z$ form a clique of four vertices in $L$,
and since $L$ is flag, they span a $3$-simplex of the $2$-dimensional
complex $L$, which is absurd. Hence $w\in A$, and the edges
$wx,wy,wz$ show that
\[
 w\in S_x\cap S_y\cap S_z.
\]
Thus the intersection is nonempty.

Conversely, let $a\in S_x\cap S_y\cap S_z$.
The vertices $a,x,y,z$ are pairwise adjacent: the three edges among
$x,y,z$ belong to the triangle, and the other three come from the
definition of $S_x,S_y,S_z$. Since $T$ is flag, $axyz$ is a
tetrahedron. Since $a\in A$, we have $a\ne v$, and since $xyz$ belongs
to only two tetrahedra, $axyz=wxyz$, i.e., $a=w$. We have proved
$S_x\cap S_y\cap S_z=\{w\}$, which is~\eqref{uni:label-triangle}.
Finally, every edge $xy$ of $L$ belongs to a triangle $xyz$ of $L$,
and the nonempty set $S_x\cap S_y\cap S_z$ is contained in
$S_x\cap S_y$. This proves~\eqref{uni:label-edge}.
\end{proof}

\begin{lemma}\label{uni:leaf}
If $a$ is a leaf of the graph $T[A]^{(1)}$, with unique neighbor $b$,
then $ab$ belongs to no square of $T$.
\end{lemma}

\begin{proof}
The hypothesis means that $N_T(a)\cap A=\{b\}$; it concerns the degree
of $a$ in $T[A]^{(1)}$, not its degree in the whole graph of $T$.
Suppose that $a-b-y-x-a$ is a square. Then its four vertices are
distinct, any two consecutive ones are joined by an edge, and the
diagonals $ay,bx$ are absent.

Since $ax$ is an edge and $a\in A$, we have $x\ne v$. Moreover
$x\notin A$, since the only neighbor of $a$ in $A$ is $b$ and $x\ne b$.
Hence $x\in V(L)$, with $a\in S_x$ and $b\notin S_x$.
By~\eqref{uni:label-connected}, the induced graph on $S_x$ is
connected. In this graph $a$ is isolated, because its only neighbor in
$A$ is $b$, which does not lie in $S_x$. Therefore
\[
 S_x=\{a\},
\]
which also follows from Lemma~\ref{lem:leaf-connected}.

Since $by$ is an edge and $b\in A$, we have $y\ne v$. If $y\in A$, the
edge $xy$ would give $y\in S_x=\{a\}$, contradicting $y\ne a$. Thus
$y\in V(L)$. Since the diagonal $ay$ is absent, $a\notin S_y$, and
consequently $S_x\cap S_y=\varnothing$. On the other hand, since
$x,y\in V(L)$ and $L$ is full, the edge $xy$ belongs to $L$.
This contradicts~\eqref{uni:label-edge}.
\end{proof}

\begin{proposition}\label{uni:degree}
The maximum degree of $T$ is eight. Moreover, $A$ consists of three
vertices $a,b,c$ spanning a triangle of $T$.
\end{proposition}

\begin{proof}
By the Davis--Okun inequality (Theorem~\ref{thm:davis-okun}; see
\cite[Theorem~11.2.1]{DavisOkun2001}),
\[
 \gamma_2(T)=16-5f_0(T)+f_1(T)\ge0.
\]
Since $f_0(T)=12$, it follows that $f_1(T)\ge44$.
Since the sum of all vertex degrees is $2f_1(T)$, the maximum degree
satisfies
\[
 8=\left\lceil\frac{88}{12}\right\rceil
 \le \left\lceil\frac{2f_1(T)}{12}\right\rceil\le d\le11.
\]
We exclude the values $11,10,9$ in turn. The cases of one or two
nonneighbors are also treated in
\cite[Lemmas~3.2(iii) and~3.4(i)]{LabbeNevo2017}.

Suppose that $d=11$. Then $v$ is adjacent to every other vertex, and
since $T$ is flag, the union of $v$ with any simplex is again a simplex.
Thus $T$ is a cone, which cannot triangulate $S^3$.

Suppose that $d=10$. Then $A=\{w\}$. By~\eqref{uni:label-connected},
every $S_z$ is a nonempty subset of $A$, so that $S_z=\{w\}$ for every
$z\in V(L)$. Hence $v$ and $w$ are nonadjacent, and both are adjacent to
every vertex of $L$. Thus the graph of $T$ is the graph of the
suspension $L*S^0_{\{v,w\}}$, and Proposition~\ref{prop:flag-clique}
gives
\[
 T=L*S^0_{\{v,w\}}.
\]
Now let $xy$ be any edge of $L$. By assumption, $xy$ lies in a square
of $T$. Each of the other two vertices of such a square is nonadjacent
to one of $x,y$. Therefore neither $v$ nor $w$ occurs in the square, since
both are adjacent to $x$ and $y$. Thus all four vertices lie in $L$,
and, $L$ being full, the square is a square of $L$. Consequently every
edge of the flag $2$-sphere $L$ belongs to a square in $L$. By
Theorem~\ref{thm:octahedral}, due to Chudnovsky and
Nevo~\cite[Lemma~3.6]{ChudnovskyNevo2020}, $L$ is the octahedral
sphere. This is impossible, since the octahedral sphere has six
vertices, whereas $f_0(L)=d=10$.

Suppose that $d=9$. Then $A$ has two vertices. By
Lemma~\ref{uni:nonneighbors} the graph induced on $A$ is connected, so
that the two vertices are joined by an edge and both are leaves of this
graph. By Lemma~\ref{uni:leaf}, this edge lies in no square of $T$,
contrary to the assumption that every edge of $T$ lies in a square.

It follows that $d=8$ and $|A|=3$. By Lemma~\ref{uni:nonneighbors} the
graph $T[A]^{(1)}$ is connected, and by Lemma~\ref{uni:leaf} and the
assumption on squares it has no leaf. A connected graph on three
vertices is either a path, which has leaves, or a triangle. Hence
$T[A]^{(1)}$ is a triangle, and since $T$ is flag, it spans a
$2$-simplex of $T$.
\end{proof}

\subsection{Three paths and the vertex count}

Let $A=\{a,b,c\}$ be as in Proposition~\ref{uni:degree}.
By Lemma~\ref{lem:sphere-links}, the link of the triangle $abc$
triangulates $S^0$, i.e., it consists of two vertices $p,q$.
Neither of them is $v$, since $v$ has no neighbors in $A$; hence they
lie in $L$. Since the link is full by Lemma~\ref{lem:link-identities},
$pq$ is not an edge of $T$. Moreover, by the
description~\eqref{eq:link-full-neighbors} of links in a flag complex
by common neighbors,
\[
 S_z=\{a,b,c\}\quad\Longleftrightarrow\quad z\in\{p,q\}.
\]

By Lemma~\ref{lem:flag-links}, every vertex of the flag $2$-sphere $L$
has degree at least four. On the other hand,
\eqref{uni:degree-decomposition} gives, for each $z\in\{p,q\}$,
\[
 \deg_T(z)=1+\deg_L(z)+3.
\]
Since the maximum degree in $T$ is eight, it follows that
$4\le\deg_L(z)\le8-4=4$ for each of the two vertices. Consequently
\begin{equation}\label{uni:pq-degree}
 \deg_L(p)=\deg_L(q)=4.
\end{equation}

We define
\[
 \begin{aligned}
 U_{ab}&=\{z\in V(L):S_z=\{a,b\}\},\\
 U_{bc}&=\{z\in V(L):S_z=\{b,c\}\},\\
 U_{ca}&=\{z\in V(L):S_z=\{c,a\}\},
 \end{aligned}
\]
and put $R=\{z\in V(L):|S_z|=1\}$.
Since every $S_z$ is a nonempty subset of the set $A$ of three
elements, the cases $|S_z|=3,2,1$ give the disjoint decomposition
\begin{equation}\label{uni:partition}
 V(L)=\{p,q\}\sqcup U_{ab}\sqcup U_{bc}\sqcup U_{ca}\sqcup R.
\end{equation}

\begin{lemma}\label{uni:paths}
For each $ij\in\{ab,bc,ca\}$, the full subcomplex
$L[\{p,q\}\cup U_{ij}]$ is an induced $p$--$q$ path whose internal
vertices are exactly $U_{ij}$. Each $U_{ij}$ is nonempty.
\end{lemma}

\begin{proof}
Consider $J=\lk_T(ab)$. By the iterated-link
identity~\eqref{eq:iterated-link},
\[
 J=\lk_{\lk_T(a)}(b).
\]
Here $\lk_T(a)$ is a triangulated $2$-sphere, and the link of its
vertex $b$ triangulates a circle. Thus $J$ is a cycle; see
Lemmas~\ref{lem:link-identities} and~\ref{lem:sphere-links}.

Since $T$ is flag, a vertex belongs to $J$ if and only if it is a
common neighbor of $a$ and $b$. The only such vertex in $A$ is $c$, and
$v$ is not such a vertex. A vertex $z\in V(L)$ is a common neighbor of
$a$ and $b$ if and only if $\{a,b\}\subseteq S_z$, i.e., $z\in U_{ab}$
or $z\in\{p,q\}$. Since $J$ is full, we obtain
\[
 J=T[\{c,p,q\}\cup U_{ab}].
\]
Moreover,
\[
 \lk_J(c)=\lk_T(abc)=\{p,q\},
\]
so that the two neighbors of $c$ in this cycle are $p$ and $q$.
We may therefore write the cycle as
\[
 c-p-u_1-\cdots-u_k-q-c,\qquad
 \{u_1,\ldots,u_k\}=U_{ab}.
\]
Removing $c$ and the edges incident to it, we obtain the path
$p-u_1-\cdots-u_k-q$. Since $J$ is full, there are no further edges of
$T$ among the remaining vertices. All of them lie in $L$, so that this
path is the full subcomplex $L[\{p,q\}\cup U_{ab}]$.
Since $p$ and $q$ are not adjacent, $k$ cannot be zero, i.e.,
$U_{ab}\ne\varnothing$. The other two cases follow by the same argument
with $a,b,c$ permuted cyclically.
\end{proof}

We denote these paths by $P_{ab},P_{bc},P_{ca}$. Their sets of
internal vertices are the pairwise disjoint sets $U_{ab},U_{bc},U_{ca}$.
Hence any two of the paths have only the vertices $p,q$ in common, and
$pq$ is not an edge. Since the paths are subcomplexes of the same
triangulated sphere, the realizations of any two of them meet exactly
in $\{p,q\}$. Thus the open arcs $|P_{ij}|\setminus\{p,q\}$ are
pairwise disjoint.

\begin{lemma}\label{uni:cross-edges}
The edge $ab$ belongs to a square if and only if there exist
$z\in U_{bc}$ and $w\in U_{ca}$ with $zw\in L^{(1)}$.
Consequently, for every two distinct sets among
$U_{ab},U_{bc},U_{ca}$, there is an edge of $L$ with one endpoint in
each set.
\end{lemma}

\begin{proof}
Suppose that $a-b-z-w-a$ is a square. Its edges and missing diagonals
give
\[
 bz,aw,zw\in T^{(1)},\qquad az,bw\notin T^{(1)}.
\]
Neither $z$ nor $w$ is $v$, because neither $a$ nor $b$ is adjacent to
$v$; and neither of them is the remaining vertex $c$ of $A$, since $c$
is adjacent to both $a$ and $b$. Since the four vertices of the square
are distinct, it follows that $z,w\in V(L)$, and
\[
 b\in S_z,\quad a\notin S_z,\qquad
 a\in S_w,\quad b\notin S_w.
\]
Since $zw$ is an edge of $L$, \eqref{uni:label-edge} gives
$S_z\cap S_w\ne\varnothing$, and the only possible common element of
these subsets of $\{a,b,c\}$ is $c$. Hence
$S_z=\{b,c\}$ and $S_w=\{a,c\}$, as required.

Conversely, suppose that $z\in U_{bc}$, $w\in U_{ca}$, and $zw$ is an
edge. Then $ab,bz,zw,wa$ are edges, while $az$ and $bw$ are not, by the
definitions of $U_{bc},U_{ca}$. Since the four vertices are distinct,
the full subcomplex on $\{a,b,z,w\}$ is the square $a-b-z-w-a$.

Since every edge of $T$ lies in a square, so does each of $ab,bc,ca$.
Applying the equivalence just proved to these three edges, we obtain,
respectively, edges between
\[
 U_{bc}\text{ and }U_{ca},\qquad
 U_{ca}\text{ and }U_{ab},\qquad
 U_{ab}\text{ and }U_{bc}.
\]
\end{proof}

The following separation lemma will be used to exclude a fourth path
with an internal vertex in $R$.

\begin{lemma}\label{uni:planar}
Let $P_1,\ldots,P_k$, with $k\ge2$, be simple $p$--$q$ paths in a
triangulation of $S^2$, and suppose that
$|P_i|\cap|P_j|=\{p,q\}$ for $i\ne j$.
Then the paths have a cyclic order such that the complement of their
union consists of $k$ open disks, each bounded by two consecutive
paths. If an embedded edge $e=zw$ joins internal vertices of two
distinct paths and
\[
 e^\circ\cap\bigcup_{i=1}^k|P_i|=\varnothing,
 \qquad e^\circ=|e|\setminus\{z,w\},
\]
then the two paths are consecutive in this cyclic order.
\end{lemma}

\begin{proof}
The initial segments of the paths meet a small circle around $p$ in a
cyclic order; here the last and the first path are also regarded as
consecutive. The union of the first two paths is a polygonal simple
closed curve, and by the Jordan--Schoenflies theorem its two
complementary regions are open disks; see
\cite[Section~2.B, p.~169]{Hatcher2002}.
The interior of each further path is connected and disjoint from the
preceding paths, and hence lies in one of the complementary disks; the
path is an arc across this disk from $p$ to $q$. Applying the same
theorem to the union of this arc with each of the two boundary paths,
we see that it divides the disk into two disks, and in the cyclic order
it lies between the two paths bounding the original disk. Both the
number of complementary disks and the description of their boundaries
now follow by induction.

To prove the last assertion, note that the connected set $e^\circ$ lies
in one complementary disk $D$. Both endpoints of $e$ are limits of
points of $e^\circ$, and hence lie on $\partial D$, which is the union
of two consecutive paths. Since each endpoint is an internal vertex of
exactly one of the paths, the paths containing $z$ and $w$ are these two
paths.
\end{proof}

\begin{proposition}\label{uni:budget}
After a permutation of $a,b,c$, one has
\[
 |U_{ab}|=2,\qquad |U_{bc}|=|U_{ca}|=1,\qquad |R|=2.
\]
Moreover, the unique neighbors of $p$ and $q$ in $R$ are distinct.
\end{proposition}

\begin{proof}
In each induced path of Lemma~\ref{uni:paths}, the endpoint $p$
has exactly one neighbor among the internal vertices.
Thus $p$ has exactly one neighbor in each of the three sets $U_{ij}$,
and these three neighbors are distinct. Since $pq$ is not an edge, in
view of~\eqref{uni:pq-degree} and the partition~\eqref{uni:partition}
the fourth neighbor of $p$ is a vertex $x\in R$, which is the unique
neighbor of $p$ in $R$. Similarly, $q$ has a unique neighbor $y\in R$.

Suppose that $x=y$. Then $p-x-q$ is a fourth simple $p$--$q$ path,
and its interior is disjoint from those of $P_{ab},P_{bc},P_{ca}$
because $x\in R$. Any three positions in a cyclic order of four contain
a nonconsecutive pair; hence two of the three original paths are
nonconsecutive in the cyclic order of the four paths. Let $P,P'$ be
such paths from $\{P_{ab},P_{bc},P_{ca}\}$.
By Lemma~\ref{uni:cross-edges} there is an edge $e=zw$ of $L$ such that
$z$ is an internal vertex of $P$ and $w$ is an internal vertex of $P'$.

We verify the disjointness condition of Lemma~\ref{uni:planar}.
Since $e$ is a $1$-simplex, its relative interior
$e^\circ=\operatorname{relint}|e|=|e|\setminus\{z,w\}$ contains no
vertex of the triangulation.
If $e^\circ$ met one of the four paths, then, since two simplices
intersect in their common face, this path would contain the whole edge
$e$, and hence both $z$ and $w$. This is impossible, since the four
paths have pairwise disjoint sets of internal vertices and $z,w$ belong
to different ones. Hence
\[
 e^\circ\cap
 \bigl(|P_{ab}|\cup|P_{bc}|\cup|P_{ca}|\cup|p-x-q|\bigr)
 =\varnothing.
\]
By Lemma~\ref{uni:planar}, $P$ and $P'$ are consecutive, contrary to
their choice. Thus $x\ne y$, and $|R|\ge2$.

Each $U_{ij}$ is nonempty. Suppose that
$|U_{ab}|+|U_{bc}|+|U_{ca}|=3$, and let $r,t,u$ be the unique vertices
of $U_{ab},U_{bc},U_{ca}$, respectively. By
Lemma~\ref{uni:cross-edges}, $rt,tu,ur$ are edges, and since $L$ is
flag, $rtu$ is a triangle. But
\[
 S_r\cap S_t\cap S_u
 =\{a,b\}\cap\{b,c\}\cap\{c,a\}=\varnothing,
\]
contrary to~\eqref{uni:label-triangle}.
Therefore $|U_{ab}|+|U_{bc}|+|U_{ca}|\ge4$.

By Proposition~\ref{uni:degree}, $f_0(L)=8$, and the
partition~\eqref{uni:partition} gives
\[
 8=2+|U_{ab}|+|U_{bc}|+|U_{ca}|+|R|
 \ge2+4+2=8.
\]
Hence equality holds throughout: $|R|=2$, and the three positive
integers $|U_{ab}|,|U_{bc}|,|U_{ca}|$ have sum four, so that they are
$2,1,1$ in some order. After a permutation of $a,b,c$ we may assume
that the set with two elements is $U_{ab}$.
\end{proof}

\subsection{Determination of the graph}

\begin{proof}[Proof of Theorem~\ref{thm:uniqueness}]
By Proposition~\ref{uni:budget}, we may write
\[
 U_{ab}=\{r,s\},\qquad U_{bc}=\{t\},\qquad
 U_{ca}=\{u\},\qquad R=\{x,y\},
\]
where $x$ is adjacent to $p$ and $y$ is adjacent to $q$.
The sets $S_z$ known so far are
\[
 S_p=S_q=\{a,b,c\},\qquad S_r=S_s=\{a,b\},\qquad
 S_t=\{b,c\},\qquad S_u=\{a,c\},
\]
while $S_x,S_y$ are singletons which remain to be determined.
We label $r,s$ along $P_{ab}$ so that the three induced paths are
\begin{equation}\label{uni:ordered-paths}
 p-r-s-q,\qquad p-t-q,\qquad p-u-q.
\end{equation}
Together with $px,qy$, these edges give all four neighbors
of $p$ and of $q$ in $L$. Thus
\begin{equation}\label{uni:pq-neighbors}
 N_L(p)=\{r,t,u,x\},\qquad N_L(q)=\{s,t,u,y\}.
\end{equation}
In particular, $pq,ps,py,qr,qx$ are not edges.
Since $U_{bc}$ and $U_{ca}$ have only one vertex each,
Lemma~\ref{uni:cross-edges} shows that $tu$ is an edge.

The link of $p$ in the flag $2$-sphere $L$ is a full cycle on the four
neighbors of $p$:
\[
 \lk_L(p)=L[\{r,t,u,x\}]\cong C_4.
\]
Since this cycle is full, it contains the edge $tu$, and the opposite
edge joins the remaining two vertices; hence $rx$ is also an edge.
By~\eqref{uni:label-edge}, $S_x\cap S_r\ne\varnothing$. Since $S_x$ is a
singleton and $S_r=\{a,b\}$, either $S_x=\{a\}$ or $S_x=\{b\}$.
If necessary, we interchange $a$ and $b$ and, at the same time, $t$
and $u$; this does not affect the paths and sets determined so far.
Hence we may assume that
\[
 S_x=\{a\}.
\]
Then $S_x\cap S_t=\varnothing$, so that $xt$ is not an edge
by~\eqref{uni:label-edge}.

A $4$-cycle containing the disjoint edges $tu,rx$ has as its other two
edges either $rt,ux$ or $ru,tx$. The latter is impossible, since $tx$
is not an edge. Hence the cycle is
\[
 \lk_L(p):\quad r-t-u-x-r,
\]
so that $rt,ux$ are edges and $ru$ is not.
By Lemma~\ref{uni:cross-edges} there is an edge between
$U_{ab}=\{r,s\}$ and $U_{ca}=\{u\}$. The only candidates are $ru,su$,
and since $ru$ is not an edge, $su$ is an edge.

Similarly, the link of $q$ is the full $4$-cycle on $\{s,t,u,y\}$.
It contains $su$ and $tu$, so that the two neighbors of $u$ in this
cycle are $s,t$, and $y$ must be adjacent to both $s$ and $t$. Thus
\[
 \lk_L(q):\quad s-u-t-y-s,
\]
so that $sy,ty$ are edges and the diagonals $st,uy$ are not.
By~\eqref{uni:label-edge}, the singleton $S_y$ meets both
$S_s=\{a,b\}$ and $S_t=\{b,c\}$. Hence its element lies in
$\{a,b\}\cap\{b,c\}=\{b\}$, that is,
\[
 S_y=\{b\}.
\]
Since $S_x\cap S_y=\varnothing$, $xy$ is not an edge.

Now $x$ is adjacent to $p,r,u$ and not adjacent to $q,t,y$. Among the
other seven vertices of $L$, only $s$ remains. Since the degree of $x$
is at least four by Lemma~\ref{lem:flag-links}, $xs$ is an edge.
Likewise, $y$ is adjacent to $q,s,t$, is not adjacent to $p,u,x$,
and has degree at least four; hence $yr$ is an edge.
The edges obtained so far are listed in
Table~\ref{uni:link-table}.

\begin{table}[tbp]
 \centering
 \renewcommand{\arraystretch}{1.12}
 \begin{tabular}{@{}ccl@{}}
 \toprule
 Vertex $z$ & $S_z=N_T(z)\cap A$ & $N_L(z)$\\
 \midrule
 $p$ & $\{a,b,c\}$ & $r,t,u,x$\\
 $q$ & $\{a,b,c\}$ & $s,t,u,y$\\
 $r$ & $\{a,b\}$ & $p,s,t,x,y$\\
 $s$ & $\{a,b\}$ & $q,r,u,x,y$\\
 $t$ & $\{b,c\}$ & $p,q,r,u,y$\\
 $u$ & $\{c,a\}$ & $p,q,s,t,x$\\
 $x$ & $\{a\}$ & $p,r,s,u$\\
 $y$ & $\{b\}$ & $q,r,s,t$\\
 \bottomrule
 \end{tabular}
 \par\medskip
 \caption{The neighbors in $A$ and in $L$ of the eight vertices of $L$.}
 \label{uni:link-table}
\end{table}

The table contains eighteen distinct edges, since the sum of its eight
degrees is $4+4+5+5+5+5+4+4=36$.
For a triangulated $2$-sphere,
$3f_2(L)=2f_1(L)$ and $f_0(L)-f_1(L)+f_2(L)=2$.
Consequently
\[
 f_1(L)=3f_0(L)-6=18.
\]
Thus the edges found above are all the edges of $L$, and the table
gives the complete neighbor sets in $L$.

The neighbor sets in $T$ are now determined as well.
For $z\in V(L)$, \eqref{uni:degree-decomposition} gives
$N_T(z)=\{v\}\sqcup N_L(z)\sqcup S_z$. For the remaining four vertices,
the triangle $abc$ and the definition of the sets $S_z$ give
\[
 \begin{aligned}
 N_T(v)&=\{p,q,r,s,t,u,x,y\},\\
 N_T(a)&=\{b,c,p,q,r,s,u,x\},\\
 N_T(b)&=\{a,c,p,q,r,s,t,y\},\\
 N_T(c)&=\{a,b,p,q,t,u\}.
 \end{aligned}
\]
These formulas determine all adjacencies and nonadjacencies in $T$.

We define a bijection $\varphi$ from $V(T)$ onto the vertex set of
$T_{12}$, numbered as in Table~\ref{tab:T12} and
\cite[Figure~12]{Milizia2024}, by
\begin{equation}\label{uni:isomorphism}
 \begin{array}{c|rrrrrrrrrrrr}
 z&v&a&b&c&p&q&r&s&t&u&x&y\\\hline
 \varphi(z)&0&5&7&10&8&4&6&3&1&2&11&9
 \end{array}.
\end{equation}
Substituting into the neighbor sets above and into
Table~\ref{uni:link-table}, we obtain
\begin{equation}\label{uni:neighbor-identification}
 \varphi\bigl(N_T(z)\bigr)=N_{T_{12}}\bigl(\varphi(z)\bigr)
 \qquad(z\in V(T)),
\end{equation}
where the sets on the right-hand side are those of
Table~\ref{tab:T12}. For example,
\[
 \varphi\bigl(N_T(x)\bigr)
 =\varphi(\{v,a,p,r,s,u\})
 =\{0,2,3,5,6,8\}=N_{T_{12}}(11),
 \qquad \varphi(x)=11.
\]
By~\eqref{uni:neighbor-identification}, $\varphi$ preserves adjacency
and nonadjacency of every pair of vertices. Since both complexes are
flag, a set of vertices spans a simplex in either complex if and only
if it is a clique in its graph. Therefore, by
Proposition~\ref{prop:flag-clique}, this graph isomorphism is an
isomorphism $T\cong T_{12}$.
\end{proof}

\begin{remark}\label{uni:numerical-consequences}
The graph obtained has three vertices of degree six and nine of degree
eight. Consequently
\[
 (f_0(T),f_1(T),f_2(T),f_3(T))=(12,45,66,33),
 \qquad \gamma_2(T)=1.
\]
These values are consequences of the classification; the proof used
only $\gamma_2(T)\ge0$, and neither its strict positivity nor a
prescribed degree sequence.
\end{remark}

\section{A reflection criterion for the boundedness of the tangent Euler class}
\label{sec:euler}

In this section we prove Theorem~\ref{thm:criterion}. We fix a flag
$3$-sphere $K$ on $I=\{0,\ldots,m-1\}$ and a nondegenerate symmetric
bilinear form $B:V\times V\to\RR$, where $V=\RR^m$, which satisfies
\eqref{eq:criterion-assumptions} and has inertia $(4,m-4)$.
We denote by $e_i$ the standard coordinate vectors, and in matrix
notation we identify $B$ with its Gram matrix
$(B(e_i,e_j))_{i,j\in I}$, so that $B(x,y)=x^{\mathsf T}By$ for
column vectors. We write $W=W_K$, $\Gamma=\Gamma_K$, $\Sigma=\Sigma_K$,
and $M=M(K)$, and denote by $p_\Sigma:\Sigma\to M$ the covering
projection of Theorem~\ref{thm:davis-manifold}.

\subsection{Bilinear forms, orthogonal complements, and projections}
Besides $B$ we use the standard Euclidean inner product and norm
\begin{equation}\label{eq:Euclidean-norm}
\langle x,y\rangle_{\mathrm E}=\sum_{i=0}^{m-1}x_i y_i,
\qquad
\|x\|_{\mathrm E}=\left(\sum_{i=0}^{m-1}x_i^2\right)^{1/2}.
\end{equation}
The form $B$ is indefinite and does not define a norm on $V$.

\begin{definition}
\label{def:bilinear-geometry}
A subspace $P\subset V$ is said to be \emph{$B$-positive} if
$B(p,p)>0$ for every $p\in P\setminus\{0\}$, and \emph{$B$-negative}
if $B(p,p)<0$ for every such $p$.
The \emph{$B$-orthogonal complement} of a subspace $P$ is
\[
P^{\perp_B}=\{n\in V:B(n,p)=0\text{ for all }p\in P\}.
\]
For a vector $v$ we write $v^{\perp_B}=(\RR v)^{\perp_B}$.
On a positive subspace, and only on such a subspace, we use the norm
$\|p\|_B=B(p,p)^{1/2}$.
The orthogonal group and the special orthogonal group of $B$ are
\[
\begin{aligned}
\OO(B)&=\{g\in\operatorname{GL}(V):B(gx,gy)=B(x,y)\text{ for all }x,y\},\\
\operatorname{SO}(B)&=\{g\in\OO(B):\det g=1\}.
\end{aligned}
\]
A vector $v$ is said to be \emph{nonisotropic} if $B(v,v)\ne0$.
\end{definition}

\begin{lemma}
\label{lem:orthogonal-projection}
Let $P\subset V$ be a $B$-positive $4$-dimensional subspace. Then
\[
V=P\oplus P^{\perp_B},\qquad B|_{P^{\perp_B}}<0.
\]
For the unique decomposition $v=p+n$ with $p\in P$ and
$n\in P^{\perp_B}$, set $\pi_P(v)=p$.
If the columns of an $m\times4$ matrix $U$ form a basis of $P$, then
\begin{equation}\label{eq:orthogonal-projection}
\pi_P(v)=U(U^{\mathsf T}BU)^{-1}U^{\mathsf T}Bv.
\end{equation}
Thus $\pi_P$ is the projection along the $B$-orthogonal complement;
in general it is not the Euclidean orthogonal projection.
\end{lemma}
\begin{proof}
The matrix $U^{\mathsf T}BU$ is positive definite, and hence
invertible. If $p$ denotes the right-hand side
of~\eqref{eq:orthogonal-projection}, then
\[
U^{\mathsf T}B(v-p)=0.
\]
Thus $p\in P$ and $v-p\in P^{\perp_B}$. Moreover,
$P\cap P^{\perp_B}=\{0\}$, since every vector in this intersection
satisfies $B(p,p)=0$. This proves the decomposition and the formula.
To show that $B$ is nondegenerate on $P^{\perp_B}$, suppose that
$n\in P^{\perp_B}$ and $B(n,n')=0$ for all $n'\in P^{\perp_B}$.
Then for every $v=p'+n'\in P\oplus P^{\perp_B}$ we have
$B(n,v)=B(n,p')+B(n,n')=0$, and the nondegeneracy of $B$ on $V$
implies $n=0$.
With respect to a basis adapted to this decomposition, the matrix of
$B$ is block diagonal. Since the positive index of $B$ is four and
$B|_P$ already has four positive eigenvalues, the complementary block
is negative definite by Sylvester's law of inertia
\cite[Theorem~4.5.8]{HornJohnson2013}.
\end{proof}

Let $g\in\OO(B)$. Then $g(P^{\perp_B})=(gP)^{\perp_B}$.
If $v=p+n$ with $p\in P$ and $n\in P^{\perp_B}$, then
$gv=gp+gn$ with $gp\in gP$ and $gn\in(gP)^{\perp_B}$, and the
uniqueness of the decomposition gives $\pi_{gP}(gv)=gp=g\pi_P(v)$.
Consequently
\begin{equation}\label{eq:projection-equivariance}
\pi_{gP}g=g\pi_P,
\qquad\text{equivalently}\qquad
\pi_Pg=g\pi_{g^{-1}P}.
\end{equation}

\subsection{The coordinate cone and its positive projections}
For $\sigma\in K$, put
\[
\CC_\sigma=\left\{\sum_{i\in\sigma}t_i e_i:t_i\ge0\right\},
\qquad
\CC(K)=\bigcup_{\sigma\in K}\CC_\sigma
=\{x\in V:x_i\ge0,\ \supp(x)\in K\}.
\]
Since each $\CC_\sigma$ is closed and $K$ is finite, $\CC(K)$ is
closed in $V$. For $x\ne0$, let $r(x)=\sum_i x_i>0$.
Then $x/r(x)\in|K|$, where $|K|$ is the realization
\eqref{eq:geometric-realization}. The map $(r,u)\mapsto ru$ therefore
identifies $\CC(K)$ with
\[
\bigl([0,\infty)\times|K|\bigr)/\bigl(\{0\}\times|K|\bigr).
\]
If $h:|K|\to S^3$ is a homeomorphism, the map
$ru\mapsto r h(u)$, with $0\mapsto0$, is a homeomorphism
$\CC(K)\to\RR^4$. In particular, the cone, including its vertex,
is a topological $4$-manifold.

\begin{lemma}\label{lem:cone-projection}
Let $P\subset V$ be a $B$-positive $4$-dimensional subspace, and
let $\pi_P$ be the projection defined in
Lemma~\ref{lem:orthogonal-projection}. Then
\[
\pi_P|_{\CC(K)}:\CC(K)\longrightarrow P
\]
is a homeomorphism.
\end{lemma}
\begin{proof}
If $x\in\CC(K)$, any two distinct indices in $\supp(x)$ form an
edge of $K$. Hence the conditions on the diagonal entries and on the
edges give
\[
B(x,x)=\sum_i B_{ii}x_i^2+
\sum_{i\ne j}B_{ij}x_i x_j=\sum_i x_i^2=\|x\|_{\mathrm E}^2.
\]
For $x,y\in\CC(K)$, all $x_i,y_j$ are nonnegative, so that
\[
B(x,y)=\sum_i x_i y_i+\sum_{i\ne j}B_{ij}x_i y_j
\le\sum_i x_i y_i=\langle x,y\rangle_{\mathrm E}.
\]
Expanding the square, we obtain
\begin{equation}\label{eq:secant-inequality}
\begin{aligned}
B(x-y,x-y)
&=B(x,x)+B(y,y)-2B(x,y)\\
&\ge\|x\|_{\mathrm E}^2+\|y\|_{\mathrm E}^2
-2\langle x,y\rangle_{\mathrm E}
=\|x-y\|_{\mathrm E}^2.
\end{aligned}
\end{equation}
Write $x-y=p+n$ according to
$V=P\oplus P^{\perp_B}$. Since $p$ and $n$ are $B$-orthogonal and $B$
is negative definite on $P^{\perp_B}$, we have
\[
B(x-y,x-y)=B(p,p)+B(n,n)\le B(p,p).
\]
Consequently,
\begin{equation}\label{eq:projection-expansion}
\|\pi_P(x-y)\|_B^2\ge B(x-y,x-y)
\ge\|x-y\|_{\mathrm E}^2.
\end{equation}
In particular, $\pi_P(x)=\pi_P(y)$ implies $x=y$.

Thus the restriction of $\pi_P$ to $\CC(K)$ is a continuous injective
map between topological $4$-manifolds. By invariance of domain
\cite[Theorem~2B.3]{Hatcher2002}, applied in charts, it is a
homeomorphism onto an open subset of $P$.
The image is also closed. In fact, suppose that $\pi_P(x_j)\to p$ with
$x_j\in\CC(K)$. Taking $y=0$ in
\eqref{eq:projection-expansion}, we see that $(x_j)$ is bounded in the
Euclidean norm. Since the cone is closed, a subsequence converges to
some $x\in\CC(K)$, and $\pi_P(x)=p$ by continuity.
Hence the image is nonempty, open, and closed in the connected
space $P$, and so it is equal to $P$.
\end{proof}

For every $g\in\OO(B)$, \eqref{eq:projection-equivariance} shows that
$\pi_P$ is also injective on $g\CC(K)$: after applying $g^{-1}$, this
is the preceding lemma for the positive subspace $g^{-1}P$.
Note that the estimate~\eqref{eq:secant-inequality} holds for arbitrary
pairs $x,y$ in the cone. This is needed below, since we shall apply it
to sets meeting several cubical faces.

\subsection{An equivariant cubical map}
For a nonisotropic vector $v$, the reflection in the $B$-orthogonal
hyperplane $v^{\perp_B}$ is
\[
r_v(u)=u-\frac{2B(u,v)}{B(v,v)}v.
\]
It maps $v$ to $-v$ and fixes $v^{\perp_B}$ pointwise; hence it
preserves $B$ and has determinant $-1$.
Since $B(e_i,e_i)=1$, we have
\begin{equation}\label{eq:reflections}
r_i(u)=r_{e_i}(u)=u-2B(u,e_i)e_i.
\end{equation}
If $ij$ is an edge of $K$, then $B(e_i,e_j)=0$, and hence
\[
r_i r_j(u)=u-2B(u,e_j)e_j-2B(u,e_i)e_i=r_j r_i(u).
\]
Together with $r_i^2=\id$, these identities show that the $r_i$ satisfy
all the defining relations of $W$. Hence we obtain a homomorphism
\[
\rho:W\longrightarrow\OO(B),\qquad\rho(s_i)=r_i.
\]
We shall not need to know whether $\rho$ is injective.

Let $\mathbf1=(1,\ldots,1)^{\mathsf T}$ and define
\begin{equation}\label{eq:base-vector}
b=-\tfrac12 B^{-1}\mathbf1.
\end{equation}
Then $B(b,e_i)=(Bb)_i=-1/2$, and hence $r_i b=b+e_i$.
For $w\in W$ and $\sigma\in K$, we denote by $Q(w,\sigma)$ the closed
cube with vertex set $wW_\sigma$. Its parametrization based at $w$ is
\[
\phi_{w,\sigma}:[0,1]^\sigma\longrightarrow Q(w,\sigma),
\qquad
\phi_{w,\sigma}(\mathbf1_A)=ws_A,
\quad s_A=\prod_{i\in A}s_i.
\]
Here $A\subseteq\sigma$ is a set of vertex indices of the simplex
$\sigma$, $(\mathbf1_A)_i=1$ for $i\in A$ and $0$ otherwise, and the
factors of $s_A$ commute. In particular, $\phi_{w,\sigma}(0)=w$.
We define an affine map on the coordinate cube by
\begin{equation}\label{eq:cubical-F}
F_{w,\sigma}(t)=\rho(w)\left(b+\sum_{i\in\sigma}t_i e_i\right).
\end{equation}
The next lemma shows that these maps are the coordinate expressions
$F\circ\phi_{w,\sigma}$ of a single map on $\Sigma$.

\begin{lemma}\label{lem:cubical-F}
There is a continuous map $F:\Sigma\to V$ such that
$F\circ\phi_{w,\sigma}=F_{w,\sigma}$ for every cube.
It satisfies $F(w)=\rho(w)b$ at every vertex $w\in W$ and
$F(hx)=\rho(h)F(x)$ for $h\in W$ and $x\in\Sigma$.
For every $z\in\Sigma$ there is an open neighborhood $U$ of $z$
such that $\pi_PF|_U$ is injective for every positive $4$-plane $P$.
\end{lemma}
\begin{proof}
For $A\subseteq\sigma$, put $R_A=\prod_{i\in A}r_i$.
Since $r_i e_j=e_j$ for distinct $i,j\in\sigma$, we have
\[
R_A b=b+\sum_{i\in A}e_i,
\qquad
R_A e_i=\begin{cases}-e_i,&i\in A,\\e_i,&i\in\sigma\setminus A.\end{cases}
\]
With respect to the base vertex $ws_A$, the same point has the
coordinates $t'$ given by $t'_i=1-t_i$ for $i\in A$ and $t'_i=t_i$
otherwise. Then
\[
\begin{aligned}
\phi_{ws_A,\sigma}(t')&=\phi_{w,\sigma}(t),\\
F_{ws_A,\sigma}(t')
&=\rho(w)R_A\left(b+\sum_{i\in\sigma}t'_i e_i\right)
=F_{w,\sigma}(t).
\end{aligned}
\]
Thus the map on a geometric cube does not depend on the choice of the
base vertex, and its value at the vertex $ws_A$ is $\rho(ws_A)b$. If
two cubes meet in a face, their affine maps agree at all vertices of
that face, and hence on the whole face. Therefore they fit together to
a map $F$ on $\Sigma$, which is continuous since the cubes form a
locally finite closed cover.

Left multiplication by $h\in W$ satisfies
$h\phi_{w,\sigma}(t)=\phi_{hw,\sigma}(t)$; the parameter $t$ is
unchanged, because both the point and the base vertex are moved by $h$.
Since $F_{hw,\sigma}(t)=\rho(h)F_{w,\sigma}(t)$, we obtain
$F(hx)=\rho(h)F(x)$.

We denote by $\cstar{w}$ the \emph{closed cubical vertex star}, i.e.,
the union of all closed cubes containing $w$. The parametrizations
based at $w$ fit together to a homeomorphism
\[
\phi_w:\CC(K)\cap[0,1]^m\longrightarrow\cstar{w},
\qquad F\circ\phi_w(t)=\rho(w)(b+t).
\]
In fact, there is one cube $Q(w,\sigma)$ for each $\sigma\in K$, and
these cubes intersect in the coordinate faces corresponding to the
intersections of the simplices. For $t,t'$ in the truncated cone, the
equality of their images under $\pi_PF\circ\phi_w$ reduces,
by~\eqref{eq:projection-equivariance}, to
$\pi_{\rho(w)^{-1}P}(t)=\pi_{\rho(w)^{-1}P}(t')$,
and Lemma~\ref{lem:cone-projection} gives $t=t'$. Thus $\pi_PF$ is
injective on the whole closed star, including pairs of points in
different cubes.

For later use we define the open sets
\begin{equation}\label{eq:open-cubical-stars}
U_w=\Sigma\setminus
\bigcup_{\substack{Q\text{ a closed cube}\\w\notin Q}}Q.
\end{equation}
By local finiteness the union is closed, and $U_w\subseteq\cstar{w}$.
The sets $U_w$ cover $\Sigma$. In fact, given $z$, let $Q_z$ be the
unique cube whose relative interior contains $z$, and choose a vertex
$w$ of $Q_z$. If a closed cube $Q$ contains $z$, then $Q\cap Q_z$ is a
face containing a relative interior point of $Q_z$, so that
$Q_z\subseteq Q$ and $w\in Q$. Consequently $z\in U_w$, and the
injectivity on $\cstar{w}$ holds on this open neighborhood of $z$.
\end{proof}

\subsection{Identification of the tangent microbundle}
The representation $\rho|_\Gamma$ defines the flat vector bundle
\begin{equation}\label{eq:flat-xi}
\xi=(\Sigma\times V)/\Gamma\longrightarrow M,
\qquad
\gamma\cdot(x,v)=(\gamma x,\rho(\gamma)v).
\end{equation}
Its projection is $[x,v]\mapsto p_\Sigma(x)$.
For $z=p_\Sigma(x)$, the map $v\mapsto[x,v]$ identifies $V$ with
$\xi_z$. If the lift $x$ is replaced by $\gamma x$, the vector
coordinates $u,v$ are replaced by $\rho(\gamma)u,\rho(\gamma)v$. Since
$\rho(\gamma)$ preserves $B$, it follows that
\[
B_\xi([x,u],[x,v])=B(u,v)
\]
is independent of the lift. Let $s:O\to\Sigma$ be a local section of
the covering. In the trivialization $(z,v)\mapsto[s(z),v]$ the form is
written as $(z,u,v)\mapsto u^{\mathsf T}Bv$, a continuous function with
constant coefficients. These local expressions agree on overlaps, and
hence $B_\xi$ is a continuous fiberwise form of inertia $(4,m-4)$.

\begin{definition}
Let $\operatorname{Gr}_4(V)$ be the Grassmannian of $4$-dimensional
linear subspaces of $V$, with its usual topology, defined by the charts
which express the subspaces transverse to a fixed complement as graphs
of linear maps. Set
\[
\PP_B=\{P\in\operatorname{Gr}_4(V):B|_P>0\}.
\]
A point of $\PP_B$ is a subspace $P$, without a choice of
orientation; positivity refers to $B|_P$.
\end{definition}

\begin{lemma}\label{lem:positive-planes}
The space $\PP_B$ is contractible, and $(P,v)\mapsto\pi_P(v)$ is
continuous. The bundle of positive planes
\[
\mathcal P(\xi)=(\Sigma\times\PP_B)/\Gamma\longrightarrow M,
\quad \gamma(x,P)=(\gamma x,\rho(\gamma)P),
\]
admits a continuous section. Such a section is equivalent to a
continuous map $\Psi:\Sigma\to\PP_B$ with
\begin{equation}\label{eq:positive-field}
\Psi(\gamma x)=\rho(\gamma)\Psi(x).
\end{equation}
It determines a positive vector subbundle of rank $4$
\begin{equation}\label{eq:positive-xi}
\xi^+=\{(x,v):x\in\Sigma,\ v\in\Psi(x)\}/\Gamma\subseteq\xi.
\end{equation}
\end{lemma}
\begin{proof}
Fix a positive $4$-plane $P_0$ and put $N_0=P_0^{\perp_B}$.
Every positive $4$-plane $P$ satisfies $P\cap N_0=\{0\}$, so that the
projection to $P_0$ restricts to an isomorphism $P\to P_0$. There is
therefore a unique linear map $A:P_0\to N_0$ with
$P=\operatorname{graph}(A)=\{p+Ap:p\in P_0\}$.
We use the norms
\[
\|p\|_+=\sqrt{B(p,p)},\qquad \|n\|_-=\sqrt{-B(n,n)},
\qquad
\|A\|_{\mathrm{op}}=\max_{\|p\|_+=1}\|Ap\|_-.
\]
By orthogonality,
\[
B(p+Ap,p+Ap)=\|p\|_+^2-\|Ap\|_-^2.
\]
Thus $B|_{\operatorname{graph}(A)}$ is positive definite if and only if
$\|Ap\|_-<1$ for every unit vector $p\in P_0$. Since the unit sphere is
compact and $p\mapsto\|Ap\|_-$ is continuous, the maximum is attained,
and therefore these strict inequalities are equivalent to
$\|A\|_{\mathrm{op}}<1$. Therefore the graph chart identifies $\PP_B$ with
\[
\mathcal B=\{A\in\operatorname{Hom}(P_0,N_0):\|A\|_{\mathrm{op}}<1\}.
\]
For $A_0,A_1\in\mathcal B$ and $0\le t\le1$, the triangle inequality
gives
\[
\|(1-t)A_0+tA_1\|_{\mathrm{op}}
\le(1-t)\|A_0\|_{\mathrm{op}}+t\|A_1\|_{\mathrm{op}}<1.
\]
Hence $\mathcal B$ is convex, and $A\mapsto(1-t)A$ contracts it to $0$.
A fixed basis $(p_j)$ of $P_0$ gives the continuous basis
$(p_j+Ap_j)$ of $\operatorname{graph}(A)$, and substituting it
in~\eqref{eq:orthogonal-projection} we obtain the continuity of
$\pi_P$.

We construct a section over the finite CW complex $M$ by induction on
skeleta, starting with arbitrary values over the vertices.
Let $\chi:D^k\to M$ be a characteristic map for a $k$-cell, with
$k\ge1$, and suppose that the section has been defined on the
$(k-1)$-skeleton. Since $D^k$ is simply connected, $\chi$ has a lift
$\widetilde\chi:D^k\to\Sigma$, which gives a trivialization
\[
D^k\times\PP_B\longrightarrow\chi^*\mathcal P(\xi),
\qquad (z,P)\longmapsto(z,[\widetilde\chi(z),P]).
\]
In this trivialization the section on the boundary is
$u\mapsto\operatorname{graph}(A(u))$, with $u\in S^{k-1}$ and
$A(u)\in\mathcal B$. Regarding $D^k$ as the cone on $S^{k-1}$, we set
\[
\widehat A(ru)=rA(u)\quad(0<r\le1),\qquad \widehat A(0)=0.
\]
Since $\|\widehat A(ru)\|_{\mathrm{op}}<r$, this extension is
continuous at the cone point. Its graph agrees with the given section
on the boundary and equals $P_0$ at the center. Boundary points
identified by $\chi$ have the same image in $\mathcal P(\xi)$,
although their coordinates with respect to different lifts need not
agree. Therefore the extension fits together with the existing section
in the CW quotient. By induction we obtain a continuous section $s$
over $M$.

For each $x\in\Sigma$ there is a unique $\Psi(x)\in\PP_B$ with
$s(p_\Sigma(x))=[x,\Psi(x)]$. The map $\Psi$ is continuous, and the
quotient relation implies~\eqref{eq:positive-field}.
The family $\{(x,v):v\in\Psi(x)\}$ is a vector subbundle of
$\Sigma\times V$, with the continuous local frames $p_j+A(x)p_j$ given
by the graph chart, and by equivariance it descends to the subbundle
\eqref{eq:positive-xi}. These frames are not necessarily equivariant,
and in general they do not descend to a global frame on $M$.
\end{proof}

We fix such a section and the resulting bundle $\xi^+$.
Next we introduce the metrics used in the comparison of microbundles.
Each cube is given its standard Euclidean metric
$d(t,t')=(\sum_i(t_i-t'_i)^2)^{1/2}$, and $\Sigma$ is given the
induced piecewise Euclidean length metric $d_\Sigma$; thus distances
are infima of lengths of paths, measured cube by cube. Here the
Euclidean structure of the cubes is used, not the indefinite form $B$;
see \cite[Section~12.2 and Section~I.3, pp.~507--508]{Davis2008}.
The cubical action of $\Gamma$ is isometric. We put
\begin{equation}\label{eq:quotient-metric}
d_M(p_\Sigma(x),p_\Sigma(y))
=\inf_{\gamma\in\Gamma}d_\Sigma(x,\gamma y).
\end{equation}
This does not depend on the choice of the lifts, and the symmetry and
the triangle inequality follow from the invariance of $d_\Sigma$. To
show that $d_M$ is positive definite, suppose that
$d_\Sigma(x,\gamma_jy)\to0$. The set
$C=\{x,y\}\cup\{\gamma_jy:j\ge1\}$ is compact. Since the action is
proper, $\{\gamma:\gamma C\cap C\ne\varnothing\}$ is finite, so that a
subsequence of $\gamma_j$ is constant, say equal to $\gamma$.
Then $x=\gamma y$. Thus distinct orbits have positive distance.

For $t>0$, write
\[
\begin{aligned}
B_\Sigma(x,t)&=\{y\in\Sigma:d_\Sigma(x,y)<t\},\\
B_M(z,t)&=\{z'\in M:d_M(z,z')<t\}.
\end{aligned}
\]
Here $B_\Sigma$ and $B_M$ denote open balls. The quotient map is
$1$-Lipschitz and satisfies
\begin{equation}\label{eq:ball-image}
p_\Sigma(B_\Sigma(x,t))=B_M(p_\Sigma(x),t).
\end{equation}
In fact, if the infimum in~\eqref{eq:quotient-metric} is less than $t$,
then some $\gamma y$ has distance less than $t$ from $x$. Thus
$p_\Sigma$ is an open map for these metrics, and $d_M$ induces the
quotient topology of the covering.

\begin{lemma}\label{lem:uniform-neighborhoods}
There exist $a,r>0$ such that
\begin{enumerate}
\item $d_\Sigma(x,\gamma x)\ge a$ for all $x\in\Sigma$ and
all $\gamma\in\Gamma\setminus\{1\}$;
\item every open ball of radius $r$ in $\Sigma$ is contained in
$U_w$ for some vertex $w$, where $U_w$ is defined in~\eqref{eq:open-cubical-stars}.
\end{enumerate}
\end{lemma}
\begin{proof}
Choose a compact set $C$ with $\Gamma C=\Sigma$, for example a
finite union of closed cubes representing all orbits of cubes.
Suppose that (i) fails. Translating and conjugating a sequence of
counterexamples, we obtain $c_j\in C$ and $h_j\ne1$ with
$d_\Sigma(c_j,h_jc_j)\to0$.
Passing to a subsequence, we may assume that $c_j\to c\in C$; then
$h_jc_j\to c$. The set $D=C\cup\{h_jc_j:j\ge1\}$ is compact, and
$h_jD\cap D\ne\varnothing$ for every $j$. Since the action is proper,
the set of the $h_j$ is finite. A constant subsequence gives $hc=c$
with $h\ne1$, which contradicts the freeness of the action.

As was shown in the proof of Lemma~\ref{lem:cubical-F}, the open sets
$U_w$ cover $\Sigma$. Cover $C$ by finitely many balls
$B_\Sigma(z_j,r_j)$ with
$B_\Sigma(z_j,2r_j)\subseteq U_{w_j}$, and put
$r=\min_j r_j>0$. If $x\in B_\Sigma(z_j,r_j)$, then
$B_\Sigma(x,r)\subseteq B_\Sigma(z_j,2r_j)$.
Thus (ii) holds for centers in $C$, and it holds for every center in
$\Sigma$, since the family of the open sets $U_w$ is invariant under
$\Gamma$ and $\Gamma$ acts isometrically.
\end{proof}

In particular, for $0<t<a/2$ we have
\[
p_\Sigma^{-1}(B_M(p_\Sigma(x),t))
=\coprod_{\gamma\in\Gamma}B_\Sigma(\gamma x,t),
\]
and each of these balls is mapped homeomorphically onto the ball in
$M$. The disjointness follows from (i), and the surjectivity
from~\eqref{eq:ball-image}. Part (ii) gives a uniform radius for the
cover $(U_w)$; the uniformity comes from the compactness of a set of
orbit representatives, although $\Sigma$ itself is noncompact.

\begin{proposition}\label{prop:tangent-identification}
The underlying topological microbundle of $\xi^+$ is isomorphic
to the tangent microbundle $\tau M$.
\end{proposition}
\begin{proof}
Let $a,r$ be as in Lemma~\ref{lem:uniform-neighborhoods}, choose
$0<\varepsilon<\min\{r,a/3\}$, and set
\[
\begin{aligned}
D_\varepsilon&=\{(x,y)\in\Sigma\times\Sigma:
d_\Sigma(x,y)<\varepsilon\},\\
U_\varepsilon&=\{(\bar x,\bar y)\in M\times M:
d_M(\bar x,\bar y)<\varepsilon\}.
\end{aligned}
\]
Let $(\bar x,\bar y)\in U_\varepsilon$, and let $x$ be a lift of
$\bar x$. By the definition of the quotient metric there is a lift
$y$ of $\bar y$ with $d_\Sigma(x,y)<\varepsilon$. It is unique:
two such lifts $y,\gamma y$ would satisfy
$d_\Sigma(y,\gamma y)<2\varepsilon<a$, which forces $\gamma=1$.
If $(x_1,y_1),(x_2,y_2)\in D_\varepsilon$ have the same image, choose
$\gamma\in\Gamma$ with $x_2=\gamma x_1$. Then $\gamma^{-1}y_2$ and
$y_1$ are lifts of the same point and both lie within $\varepsilon$ of
 $x_1$. The uniqueness just proved gives $\gamma^{-1}y_2=y_1$.
Thus $(p_\Sigma\times p_\Sigma)|_{D_\varepsilon}$ identifies exactly
the points in the same orbit of the diagonal action
$\gamma(x,y)=(\gamma x,\gamma y)$.
This restriction is open, since $D_\varepsilon$ is open and
$p_\Sigma\times p_\Sigma$ is a covering map. Hence it induces a
homeomorphism
\[
D_\varepsilon/\Gamma\ \cong\ U_\varepsilon.
\]

Consider the pullback of $\xi^+$ along the covering map $p_\Sigma$.
Its definition and its concrete model are
\[
\begin{aligned}
p_\Sigma^*\xi^+
&=\{(x,\eta):x\in\Sigma,\ \eta\in\xi^+_{p_\Sigma(x)}\}\\
&\cong\widetilde E^+:=\{(x,v):v\in\Psi(x)\},
\qquad (x,v)\longmapsto(x,[x,v]).
\end{aligned}
\]
We write $\widetilde p(x,v)=x$ and $q_\Sigma(x,v)=[x,v]$, and denote by
$p_+:\xi^+\to M$ the bundle projection. The pullback square is
\[
\begin{tikzcd}[column sep=large]
\widetilde E^+ \arrow[r,"q_\Sigma"] \arrow[d,"\widetilde p"']
 & \xi^+ \arrow[d,"p_+"]\\
\Sigma \arrow[r,"p_\Sigma"'] & M.
\end{tikzcd}
\]
We define $\Theta:D_\varepsilon\to\widetilde E^+$ by
\begin{equation}\label{eq:Theta}
\Theta(x,y)=\bigl(x,\pi_{\Psi(x)}(F(y)-F(x))\bigr),
\qquad (x,y)\in D_\varepsilon.
\end{equation}
It is continuous by Lemma~\ref{lem:positive-planes}.
For $\gamma\in\Gamma$, \eqref{eq:projection-equivariance}
and the equivariance of $F$ and $\Psi$ give
\[
\pi_{\Psi(\gamma x)}\bigl(F(\gamma y)-F(\gamma x)\bigr)
=\rho(\gamma)\pi_{\Psi(x)}\bigl(F(y)-F(x)\bigr).
\]
Thus $\Theta$ induces a continuous map
\[
\overline\Theta:U_\varepsilon\longrightarrow\xi^+,
\qquad
(p_\Sigma(x),p_\Sigma(y))\longmapsto
[x,\pi_{\Psi(x)}(F(y)-F(x))].
\]
In other words, the following diagram commutes:
\[
\begin{tikzcd}[column sep=large]
D_\varepsilon \arrow[r,"\Theta"]
 \arrow[d,"(p_\Sigma\times p_\Sigma)|_{D_\varepsilon}"']
 & \widetilde E^+ \arrow[d,"q_\Sigma"]\\
U_\varepsilon \arrow[r,"\overline\Theta"'] & \xi^+.
\end{tikzcd}
\]
Moreover $\widetilde p\Theta=\operatorname{pr}_1$, so that
$p_+\overline\Theta=\operatorname{pr}_1$. Since
$\Theta(x,x)=(x,0)$, the diagonal is mapped to the zero section.

To prove that $\overline\Theta$ is injective, suppose that two points
of $U_\varepsilon$ have the same image. Their base points are equal; we
choose a common lift $x$ and the unique lifts $y_1,y_2$ of their
second coordinates with $d_\Sigma(x,y_i)<\varepsilon$. With this lift
fixed, the equality in the fiber means
\[
\pi_{\Psi(x)}F(y_1)=\pi_{\Psi(x)}F(y_2).
\]
Since $\varepsilon<r$, the points $x,y_1,y_2$ lie in one closed vertex
star, and the injectivity on this star proved in
Lemma~\ref{lem:cubical-F} gives $y_1=y_2$.

Both $U_\varepsilon$, an open subset of $M\times M$, and the
total space of the bundle $\xi^+\to M$ of rank $4$ are topological
$8$-manifolds. By invariance of domain
\cite[Theorem~2B.3]{Hatcher2002}, $\overline\Theta$ is a
homeomorphism onto an open subset of $\xi^+$, and this subset contains
the zero section because $\Theta(x,x)=(x,0)$.
Thus $\overline\Theta$ maps a neighborhood of the diagonal
homeomorphically onto a neighborhood of the zero section, commutes with
the projections to $M$, and preserves the sections; these are exactly
the conditions for an isomorphism of microbundles. Denoting by $0_+$
the zero section, we have the commutative diagram
\[
\begin{tikzcd}[column sep=large]
M \arrow[r,"\Delta"] \arrow[d,equal]
 & U_\varepsilon \arrow[r,"\operatorname{pr}_1"] \arrow[d,"\overline\Theta"]
 & M \arrow[d,equal]\\
M \arrow[r,"0_+"'] & \xi^+ \arrow[r,"p_+"'] & M.
\end{tikzcd}
\]
\end{proof}

Proposition~\ref{prop:tangent-identification} concerns only
neighborhoods of the sections; in particular, it requires neither the
global injectivity of $F$ nor the flatness of $\xi^+$. We next use the
connected structure group of $\xi$ to orient $\xi^+$ and construct
its universal Euler class.

\subsection{The connected structure group}\label{subsec:structure-group}
The representation $\rho|_\Gamma$ takes values in the identity
component of the special orthogonal group. This provides the connected
structure group needed to orient the positive subbundle and apply
Theorem~\ref{thm:primary-bounded}.
Put
\begin{equation}\label{eq:structure-group}
G=\operatorname{SO}_0(B),
\end{equation}
where the subscript denotes the connected component containing the
identity in the usual matrix topology.

\begin{lemma}\label{lem:identity-component}
The reflection representation satisfies $\rho(\Gamma)\subseteq G$.
\end{lemma}
\begin{proof}
The positive unit quadric
\[
Q_+=\{v\in V:B(v,v)=1\}
\]
is path connected. Indeed, after choosing coordinates in which
$B=\operatorname{diag}(I_4,-I_{m-4})$, the map
\[
S^3\times\RR^{m-4}\longrightarrow Q_+,
\qquad (u,z)\longmapsto\bigl(\sqrt{1+\|z\|^2}\,u,z\bigr)
\]
is a homeomorphism. Each original basis vector satisfies
$B(e_i,e_i)=1$, so for each $i$ we may choose a path
$v_i:[0,1]\to Q_+$ from $e_i$ to $e_0$.
The reflections $r_{v_i(t)}$ depend continuously on $t$.

Let $\gamma\in\Gamma=\ker\varepsilon$. Every generator occurs an
even number of times in a word representing $\gamma$, so its length
is even. Write $\gamma=s_{i_1}\cdots s_{i_{2k}}$. Then
\[
t\longmapsto r_{v_{i_1}(t)}\cdots r_{v_{i_{2k}}(t)}
\]
is a path in $\operatorname{SO}(B)$ from $\rho(\gamma)$ to
$r_{e_0}^{2k}=1$: each factor preserves $B$ and has determinant $-1$.
Thus $\rho(\gamma)$ belongs to $G$.
\end{proof}

A vertex link in $K$ has at least six vertices by
Lemma~\ref{lem:flag-links}. Every vertex also has a nonneighbor,
since otherwise flagness would make $K$ a cone. Hence $m\ge8$.
The group $G$ is isomorphic to $\operatorname{SO}_0(4,m-4)$, whose
Lie algebra is simple by~\cite[Theorem~6.105(c), p.~421]{Knapp2002}.
Thus $G$ is a connected semisimple Lie group, and
Theorem~\ref{thm:primary-bounded} applies to it.

\subsection{The universal positive Euler class}
We construct a class in $H^4(BG;\RR)$ whose pullback is $e(\xi^+)$.
We first orient the tautological bundle of positive planes and check
that the action of $G$ preserves this orientation.

The \emph{tautological bundle} over $\PP_B$ is
\[
 \mathcal T=\{(P,v)\in\PP_B\times V:v\in P\}
 \longrightarrow\PP_B,\qquad (P,v)\longmapsto P.
\]
Its fiber $\mathcal T_P$ is the positive $4$-plane $P$ itself. We
orient $\mathcal T$ as follows. Fix the decomposition $V=P_0\oplus N_0$
used in Lemma~\ref{lem:positive-planes} and an ordered basis
$p_1,\ldots,p_4$ of $P_0$, and write $P=\operatorname{graph}(A_P)$ as
in the graph chart. The vectors
\[
 \beta_j(P)=p_j+A_Pp_j\quad(1\le j\le4)
\]
form a basis of $\mathcal T_P$ which depends continuously on $P$. We
orient each fiber by this ordered basis; since the basis varies
continuously with $P$, this gives an orientation of $\mathcal T$.
The group $G$ acts on $\mathcal T$ by
\[
 g\cdot(P,v)=(gP,gv).
\]
Preservation of orientation means that $g:P\to gP$ preserves the
chosen orientations of these two fibers.

\begin{lemma}\label{lem:positive-orientation}
Every element of $G$ preserves the chosen orientation of $\mathcal T$.
\end{lemma}
\begin{proof}
For $g\in G$ and $P\in\PP_B$, define the matrix
$C(g,P)=(c_{ji}(g,P))$ by
\[
g\beta_i(P)=\sum_{j=1}^4c_{ji}(g,P)\beta_j(gP).
\]
For fixed $P$, this is a continuous invertible matrix as a function
of $g$, and $C(1,P)=I_4$. Since $G$ is connected,
$\det C(g,P)>0$ for every $g\in G$. Thus $g:P\to gP$ preserves
the chosen orientations.
\end{proof}

Let $EG\to BG$ be a universal principal right $G$-bundle with $EG$
contractible; here $G$ has its usual topology. By universality, every
principal $G$-bundle over the finite CW complex $M$ is isomorphic
to a pullback of this bundle; see the conventions on classifying spaces
in Section~\ref{sec:preliminaries}. The natural action of $G$ on $V$
gives the associated vector bundle of rank $m$
\[
 \mathcal E=EG\times_G V\longrightarrow BG,
 \qquad (ug,v)\sim(u,gv).
\]
Since $G$ preserves $B$, the fibers of $\mathcal E$ inherit this
form. The bundle of positive $4$-planes of $\mathcal E$ is
\[
 Y=EG\times_G\PP_B,\qquad q:Y\longrightarrow BG,
 \qquad q([u,P])=[u].
\]
The equivalence relation is $(ug,P)\sim(u,gP)$.
Intrinsically, a point of $Y$ is a pair $(b,L)$ of a point $b\in BG$
and a positive $4$-dimensional subspace $L\subseteq\mathcal E_b$; the
point $[u,P]$ corresponds to
\[
 b=[u],\qquad L=\{[u,v]:v\in P\}\subseteq\mathcal E_b.
\]
The map $q$ forgets $L$, and its fiber is $\PP_B$.

The associated tautological bundle is
\[
 \mathcal U^+=EG\times_G\mathcal T\longrightarrow Y,
 \qquad [u,(P,v)]\longmapsto[u,P].
\]
The equivalence relation in its total space is
$(ug,P,v)\sim(u,gP,gv)$. In the intrinsic description,
\[
 \mathcal U^+_{(b,L)}=L\subseteq\mathcal E_b,
 \qquad \mathcal U^+\subseteq q^*\mathcal E.
\]
Thus $\mathcal U^+$ has rank four and is defined over $Y$, whose
points include the choice of a positive $4$-plane.
By Lemma~\ref{lem:positive-orientation}, the orientation of
$\mathcal T$ descends to $\mathcal U^+$, and we obtain
$e(\mathcal U^+)\in H^4(Y;\RR)$. The frame $(\beta_j)$ trivializes
$\mathcal T$, but it is not necessarily $G$-equivariant, and in general
it does not descend to a frame of $\mathcal U^+$.

Since the fiber $\PP_B$ of $q$ is contractible, the homotopy exact
sequence shows that $q$ is a weak homotopy equivalence; see
\cite[Theorem~4.41 and Proposition~4.48]{Hatcher2002}. A weak homotopy
equivalence induces an isomorphism of singular cohomology
\cite[Proposition~4.21]{Hatcher2002}. In particular,
\[
 q^*:H^4(BG;\RR)\xrightarrow{\ \cong\ }H^4(Y;\RR).
\]

\begin{definition}
The chosen orientation determines the unique class
\begin{equation}\label{eq:universal-positive-Euler}
 \mathfrak e^+=(q^*)^{-1}e(\mathcal U^+)\in H^4(BG;\RR),
 \qquad q^*\mathfrak e^+=e(\mathcal U^+).
\end{equation}
\end{definition}
Here $(q^*)^{-1}$ is the inverse of the isomorphism $q^*$ of cohomology
groups; no map inverse to $q$ is involved.

By Lemma~\ref{lem:identity-component}, the representation
$\rho|_\Gamma:\Gamma\to G$ defines the flat principal right
$G$-bundle
\[
 \mathcal G=(\Sigma\times G)/\Gamma\longrightarrow M,
 \qquad \gamma\cdot(x,a)=(\gamma x,\rho(\gamma)a),
 \qquad [x,a]\cdot g=[x,ag].
\]
Let $c_\delta:M\to BG^\delta$ be a classifying map of $\mathcal G$ with
its discrete structure group, and put $f=\iota\circ c_\delta$, where
$\iota:BG^\delta\to BG$ is induced by the continuous identity
$G^\delta\to G$. We choose an isomorphism $\mathcal G\cong f^*EG$, or,
equivalently, a map of principal bundles for which the following square
is a pullback:
\[
\begin{tikzcd}[column sep=large]
 \mathcal G \arrow[r] \arrow[d] & EG \arrow[d]\\
 M \arrow[r,"f"'] & BG.
\end{tikzcd}
\]
The associated bundle of rank $m$ is identified with $\xi$ by
\[
 \mathcal G\times_G V\longrightarrow\xi,
 \qquad [[x,a],v]\longmapsto[x,av].
\]
Consequently the same isomorphism of principal bundles gives
$\xi\cong f^*\mathcal E$, and identifies the bundle of positive planes
of $\xi$ with $f^*Y$. More explicitly, we have the pullback square
\[
\begin{tikzcd}[column sep=large]
 \mathcal P(\xi)\cong f^*Y \arrow[r,"\widehat f"] \arrow[d]
 &Y \arrow[d,"q"]\\
 M \arrow[r,"f"'] &BG.
\end{tikzcd}
\]
Here $\mathcal P(\xi)=\mathcal G\times_G\PP_B$, and its fiber over
$z\in M$ consists of the positive $4$-planes in $\xi_z$.

The map $\Psi:\Sigma\to\PP_B$ gives a section
\[
 s:M\longrightarrow\mathcal P(\xi),
 \qquad s(p_\Sigma(x))=[[x,1],\Psi(x)].
\]
This formula does not depend on the lift $x$. In fact,
$[\gamma x,1]=[x,\rho(\gamma)^{-1}]$ in $\mathcal G$ and
$\Psi(\gamma x)=\rho(\gamma)\Psi(x)$, and these two changes cancel by
the relation in the associated bundle. In other words, $s(z)$ is the
subspace $\xi_z^+\subseteq\xi_z$.
Put $\widetilde f=\widehat f\circ s$. Since $s$ is a section,
$q\circ\widetilde f=f$. Conversely, every such lift of $f$ defines a
section of $f^*Y$ by $z\mapsto(z,\widetilde f(z))$.

By the definition of $\mathcal U^+$, its pullback by $\widetilde f$
has fiber $\xi_z^+$ over $z$. The resulting oriented bundle isomorphism
and the maps of base spaces fit into the diagram
\[
\begin{tikzcd}[column sep=large,row sep=large]
 \xi^+ \arrow[r] \arrow[d] &\mathcal U^+ \arrow[d]\\
 M \arrow[r,"\widetilde f"] \arrow[d,equal]
 &Y \arrow[d,"q"]\\
 M \arrow[r,"f"'] &BG.
\end{tikzcd}
\]
The upper square is a pullback, and the lower square commutes. Thus
$\xi^+\cong\widetilde f^*\mathcal U^+$, with the orientation induced
from $\mathcal T$, and the naturality of the Euler class yields
\begin{equation}\label{eq:pullback-positive-Euler}
\begin{split}
 e(\xi^+)
 &=\widetilde f^*e(\mathcal U^+)
 =\widetilde f^*q^*\mathfrak e^+\\
 &=(q\circ\widetilde f)^*\mathfrak e^+
 =f^*\mathfrak e^+
 =c_\delta^*\iota^*\mathfrak e^+.
\end{split}
\end{equation}
For the fixed orientation convention, the right-hand side depends only
on the principal bundle, so that $e(\xi^+)$ does not depend on the
choice of the section of positive planes. This identity expresses
$e(\xi^+)$ as a primary characteristic class of the flat principal
$G$-bundle $\mathcal G$. It does not mean that $\xi^+$ is flat, since
parallel transport in $\xi$ need not preserve the varying planes
$\Psi(x)$; the flatness required in Theorem~\ref{thm:primary-bounded}
is that of $\mathcal G$.

\begin{proof}[Proof of Theorem~\ref{thm:criterion}]
The group $G=\operatorname{SO}_0(B)$ is connected and semisimple.
Applying Theorem~\ref{thm:primary-bounded} to
$\mathfrak e^+\in H^4(BG;\RR)$, we obtain a
bounded singular cocycle $\zeta$ on $BG^\delta$ representing
$\iota^*\mathfrak e^+$. Its pullback
$\eta=c_\delta^*\zeta$ is bounded, with
$\|\eta\|_\infty\le\|\zeta\|_\infty$, and its cohomology class is
$e(\xi^+)$ by~\eqref{eq:pullback-positive-Euler}.

By Lemma~\ref{lem:positive-orientation}, $\xi^+$ is oriented. We
transport this orientation by the isomorphism of
Proposition~\ref{prop:tangent-identification}, and thus orient the
tangent microbundle and $M$. By the naturality of the Thom class,
as in Proposition~\ref{prop:microbundle-euler}, we have
$e(\xi^+)=e(\tau M)$, and~\eqref{eq:Euler-pairing} gives
$\langle[\eta],[M]\rangle=\chi(M)$.
If $\chi(M)\ne0$, then $0<\|\eta\|_\infty<\infty$, and the
duality estimate~\eqref{eq:duality} yields
\[
\volsimp M\ge\frac{|\chi(M)|}{\|\eta\|_\infty}>0.
\]
\end{proof}

\section{The explicit form for \texorpdfstring{$T_{12}$}{T12}}\label{sec:application}

In this section we apply Theorem~\ref{thm:criterion} to $K=T_{12}$ and
prove Corollary~\ref{thm:positive}. Let $e_0,\ldots,e_{11}$ be the
standard basis of $\RR^{12}$. By means of Table~\ref{tab:T12} we
define a symmetric bilinear form by
\begin{equation}\label{eq:explicit-B}
B(e_i,e_j)=
\begin{cases}
1,&i=j,\\
0,&\{i,j\}\text{ is an edge of }K,\\
-4,&i\ne j,\ \{i,j\}\text{ is not an edge of }K.
\end{cases}
\end{equation}
We also write $B=(B_{ij})_{0\le i,j\le11}$ for its Gram matrix
$B_{ij}=B(e_i,e_j)$ in this ordered basis, so that
$B(x,y)=x^{\mathsf T}By$ for column vectors $x,y$.
Obviously $B$ satisfies the diagonal, zero, and sign conditions
in~\eqref{eq:criterion-assumptions}. Nondegeneracy and the inertia are
verified by exact integer computations.

\begin{proposition}\label{prop:explicit-inertia}
The form~\eqref{eq:explicit-B} is nondegenerate, has inertia
$(4,8)$, and has determinant $152881$.
\end{proposition}

\begin{proof}
Order the vertices as $S=(0,1,2,4)$ followed by
$J=(3,5,6,7,8,9,10,11)$. The vertices in $S$ are pairwise adjacent,
so the matrix in this basis is
\[
\widetilde B=\begin{pmatrix}I_4&C\\C^{\mathsf T}&D\end{pmatrix},
\qquad
C=(B_{sj})_{s\in S,j\in J},\quad
D=(B_{jk})_{j,k\in J},
\]
where $I_r$ denotes the identity matrix of size $r$.
Set $H=C^{\mathsf T}C-D$. Direct calculation
from~\eqref{eq:explicit-B} gives, in the specified order on $J$,
\begin{equation}\label{eq:explicit-H}
H=
\begin{pmatrix}
15&16&0&0&4&0&4&16\\
16&31&0&16&0&4&16&16\\
0&0&31&16&16&16&4&16\\
0&16&16&31&0&16&16&4\\
4&0&16&0&15&4&0&16\\
0&4&16&16&4&15&4&4\\
4&16&4&16&0&4&15&4\\
16&16&16&4&16&4&4&31
\end{pmatrix}.
\end{equation}
Let $\Delta_k$ be the determinant of the submatrix formed by the first
$k$ rows and columns of $H$. Exact computation gives
\begin{equation}\label{eq:principal-minors}
(\Delta_1,\ldots,\Delta_8)=
(15,\ 209,\ 6479,\ 28305,\ 50575,\ 157505,\ 378879,\ 152881).
\end{equation}
All these leading principal minors are positive, so $H$ is positive
definite by Sylvester's criterion
\cite[Theorem~7.2.5(b)]{HornJohnson2013}.

The matrix $R=\begin{pmatrix}I_4&-C\\0&I_8\end{pmatrix}$ has
determinant one and satisfies
\begin{equation}\label{eq:explicit-congruence}
R^{\mathsf T}\widetilde B R
=\begin{pmatrix}I_4&0\\0&-H\end{pmatrix}.
\end{equation}
Sylvester's law of inertia~\cite[Theorem~4.5.8]{HornJohnson2013}
therefore gives inertia $(4,8)$ and nondegeneracy. Since simultaneous
permutation of rows and columns preserves the determinant, we also have
\[
\det B=\det\widetilde B=\det(-H)=(-1)^8\Delta_8=152881.
\]
\end{proof}

The base vector~\eqref{eq:base-vector} can also be computed
directly. In the original vertex order $0,1,\ldots,11$, put
\[
 b=
 (\underbrace{\tfrac32,\ldots,\tfrac32}_{9\text{ entries}},
 -\tfrac52,-\tfrac52,-\tfrac52)^{\mathsf T}.
\]
By Table~\ref{tab:T12}, the numbers of nonneighbors of a vertex $i$,
the vertex itself being excluded, are as follows:
\[
 \begin{array}{c|cc}
 i&
 \text{in }\{0,\ldots,8\}&
 \text{in }\{9,10,11\}\\ \hline
 0\le i\le8&2&1\\
 9\le i\le11&3&2
 \end{array}.
\]
Since the off-diagonal entries of $B$ vanish on edges and are equal to
$-4$ on nonedges, it follows that
\[
 (Bb)_i=
 \begin{cases}
 \dfrac32-4\left(2\cdot\dfrac32-\dfrac52\right)
       =-\dfrac12,&0\le i\le8,\\[6pt]
 -\dfrac52-4\left(3\cdot\dfrac32-2\cdot\dfrac52\right)
       =-\dfrac12,&9\le i\le11.
 \end{cases}
\]
Thus $Bb=-\tfrac12\mathbf1$, where $\mathbf1$ is the column vector
with all twelve entries equal to one.
Since $B$ is nondegenerate, this $b$ coincides with
$-\tfrac12B^{-1}\mathbf1$.
Moreover, by the reflection formula~\eqref{eq:reflections}, every
$r_i$, and hence every $\rho(w)$, is given by an integer matrix.
Consequently the affine maps defining $F$
in~\eqref{eq:cubical-F} have rational coefficients.

\begin{proof}[Proof of Corollary~\ref{thm:positive}]
Put $M=M(T_{12})$.
By~\eqref{eq:explicit-B} and Proposition~\ref{prop:explicit-inertia},
all the hypotheses of Theorem~\ref{thm:criterion} are satisfied.
Therefore the Euler class $e(\tau M)\in H^4(M;\RR)$ is represented by
a bounded singular cocycle.

We compute its value on $[M]$.
By~\eqref{eq:T12-fvector}, $\gamma_2(T_{12})=45-5\cdot12+16=1$, and
hence~\eqref{eq:euler-Davis} gives
\[
 \chi(M)=2^{12-4}\gamma_2(T_{12})=256.
\]
This can also be seen directly from the cubical decomposition.
An open $r$-cell has exactly $r$ varying coordinates, whose
indices form an $(r-1)$-simplex of $T_{12}$, and each of the other
$12-r$ coordinates is fixed, with an independent choice of sign.
Thus the number $c_r$ of $r$-cells is
\[
 c_r=2^{12-r}f_{r-1}(T_{12})\qquad(0\le r\le4),
\]
where $f_{-1}=1$ corresponds to the empty simplex. Explicitly,
\[
 \begin{array}{c|rrrrr}
 r&0&1&2&3&4\\ \hline
 f_{r-1}(T_{12})&1&12&45&66&33\\
 c_r&4096&24576&46080&33792&8448
 \end{array},
\]
and their alternating sum is
\[
 4096-24576+46080-33792+8448=256.
\]
Therefore the Euler pairing~\eqref{eq:Euler-pairing} gives
\[
 \langle e(\tau M),[M]\rangle=\chi(M)=256.
\]

Let $\eta\in C_b^4(M;\RR)$ be a bounded singular cocycle whose
cohomology class is $e(\tau M)$. Its norm
$\|\eta\|_\infty=\sup_{\sigma:\Delta^4\to M}|\eta(\sigma)|$
is finite. For any real singular fundamental cycle
$z=\sum_j a_j\sigma_j$, since $\eta$ is a cocycle and $[z]=[M]$, we
have
\[
 256
 =|\eta(z)|
 =\left|\sum_j a_j\eta(\sigma_j)\right|
 \le\|\eta\|_\infty\sum_j|a_j|
 =\|\eta\|_\infty\|z\|_1.
\]
In particular, $\|\eta\|_\infty>0$. Taking the infimum over all
fundamental cycles $z$, as in~\eqref{eq:duality}, we obtain
\[
 \volsimp{M(T_{12})}
 \ge \frac{256}{\|\eta\|_\infty}>0.
\]
\end{proof}

\noindent\textbf{Acknowledgments.}
Xueyuan Wan is supported by the National Key R\&D Program of China (Grant No.~2024YFA1013200) and the National Natural Science Foundation of China (Grant No.~12671100).

\noindent\textbf{AI Declaration}: The authors used ChatGPT as an auxiliary tool in this work. The authors verified and completed all mathematical arguments and take full responsibility for the content of this paper.

\providecommand{\bysame}{\leavevmode\hbox to3em{\hrulefill}\thinspace}
\providecommand{\MR}{\relax\ifhmode\unskip\space\fi MR }
\providecommand{\MRhref}[2]{%
  \href{http://www.ams.org/mathscinet-getitem?mr=#1}{#2}
}
\providecommand{\href}[2]{#2}

\end{document}